%% file: main.tex
\documentclass[final,authoryear,jmlmc]{beg_39}             %  final version
\input{macros.tex}
\usepackage[hang]{footmisc}
\fancypagestyle{plain}{%
	\fancyhf{}
	\fancyhead[R]{\small {\it \jname}, x(x):\thepage--\pageref{LastPage} (\myyear\today)}
	\fancyfoot[R]{\small\bf\thepage }
	\fancyfoot[L]{\fottitle}
}
\usepackage{graphicx}
\usepackage{amssymb}
\usepackage{amsmath,verbatim,amsfonts,amscd}
\usepackage{amsthm}
\usepackage{dsfont}
\usepackage{hyperref}
\usepackage{latexsym,bm,color}

\renewcommand{\dmyy}{20}
\renewcommand{\myyear}{2023}
\renewcommand{\today}{}

\renewcommand{\P}{\mathbf{P}}

\theoremstyle{definition}

\begin{document}
	\volume{Volume x, Issue x, \myyear\today}
	\title{Error Estimates of Residual Minimization
		using Neural Networks 
		for Linear PDEs}
	\titlehead{Error Estimates of Residual Minimization using NNs for Linear PDEs}
	\authorhead{Y. Shin, Z. Zhang, G. E. Karniadakis}
	%For at least  authors with different addresses, use instead the following commands
	\corrauthor[1]{Yeonjong Shin}
	\author[2]{Zhongqiang Zhang}
	\author[3]{George Em Karniadakis}
	\corremail{yeonjong\_shin@ncsu.edu}
%	\corraddress{Department of Mathematics, North Carolina State University, Raleigh, NC 27695, USA}
	\address[1]{Department of Mathematics, North Carolina State University, Raleigh, NC 27695, USA}
	\address[2]{Department of Mathematical Sciences,	
		Worcester Polytechnic Institute, Worcester , MA 01609, USA}
	\address[3]{Division of Applied Mathematics,
		Brown University,
		Providence, RI 02912, USA}
	%\address[2]{Business or Academic Affiliation 2, City, Province, Zip Code, Country Business or Academic Affiliation 2, City, Province, Zip Code, Country}
	% End information for at least  authors with different addresses
	% For authors with the same post address,
	%\corrauthor{First A. Author}
	%\corremail{f.author@affiliation.com}
	%\author{Second B. Author, Jr.}
	%\address{Department of Chemistry and Courant, Institute of Mathematical Sciences, New York, NY 10012, USA}
	% End commands for all authors with the same address
	
	\dataO{mm/dd/yyyy}
	%\dataO{}
	\dataF{mm/dd/yyyy}
	%\dataF{}
	
	\abstract{
		We propose an abstract framework for analyzing the convergence of least-squares methods based on residual minimization when feasible solutions are neural networks. With the norm relations and compactness arguments, we derive error estimates for both continuous and discrete formulations of residual minimization in strong and weak forms. 
		The formulations cover recently developed physics-informed neural networks based on strong and variational formulations. 
	}
	

	\maketitle

	\input intro   %\newpage
	\input setup  %  \newpage
	\input crm   %\newpage 
	\input drm   %\newpage
	\input hpvpinn  %\newpage 
	\input discussion %.tex
	%\input inverse.tex

	%% The Acknowledgements part is started with the command \acknowledgements;
	%% acknowledgements are then done as normal sections before appendix
	%% \acknowledgements

	\acknowledgements
	This work was supported by 
	the PhiLMs grant (DE-SC0019453),
	the DARPA CompMods grant (HR00112090062), 
	the AFOSR grant (FA9550-20-1-0056),
	the MURI AFOSR grant (FA9550-20-1-0358),
	and the MURI ARO grant (W911NF-15-1-0562).

	%% appendix sections are then done as normal sections and after Acknowledgements
	
	\appendix
	\input appendix

	\input elliptic  %\newpage
	\input advection  %\newpage
	
	\bibliographystyle{Bibliography_Style}
	\bibliography{ref}
\end{document}

%% file: macros.tex
\usepackage{amsfonts}
\usepackage{amsmath,amssymb,amsthm} 

\usepackage{enumerate} 
\usepackage{graphicx}
\usepackage{xcolor}
\usepackage[hidelinks]{hyperref}
\usepackage{bm}
\let\hide\iffalse
 
\newcommand{\abs}[1]{\left\vert#1\right\vert}

\newcommand{\set}[1]{\left\{#1\right\}}
\newcommand{\Real}{\mathbb R}

\newcommand{\bN}{\mathbb{N}}

\newcommand{\bP}{\mathbb{P}}

%% file: intro.tex
\section{Introduction} \label{sec:intro}

Deep learning algorithms using neural networks 
have been employed to solve forward and inverse problems for partial differential equations (PDEs).
Many works have shown its effectiveness in various applications, such as in 
\citep{BerNys18,EYu18,khoo2019solving,Lagaris_98_ANN-ODE-PDE,Lagaris_00_ANN-Irregular,LiTWL19,liaoMing19,mao2020physics,NabMei19,pang2019fpinns,raissi2018hidden,raissi2017machine,RaiPK19,RaiYK20,SamANe20,SirSpi18,zhang2019quantifying}.
% including fractional PDEs, stochastic differential equations, biomedical problems, and fluid mechanics.
Empirical studies  \citep{raissi2018hidden,RaiYK20}
have found that deep learning algorithms are particularly effective 
in solving inverse problems for PDEs with a few data points.
Such inverse problems are known to be challenging 
for existing classical methods.
However, the mathematical theory of deep learning algorithms for PDEs is  far from being complete at the moment.

Due to the nonlinear and compositional nature of neural networks, deep learning for PDEs is often recast as highly nonconvex and nonlinear optimization problems. 
A critical step in deep learning is the formulation of an appropriate loss functional. 
A least-squares type loss is a common choice when it comes to regression or supervised learning.
%For supervised learning and regression problems,
%a least-squares type loss functional is a common choice, as in least-squares finite element methods 
The least-squares finite element methods (LS-FEM) \citep{BocGun16,BraSch70} are classic examples of using such loss functionals. 
The LS-FEM uses a linear finite element space for feasible solutions and solves linear systems for linear problems. In contrast, deep learning methods  use neural networks as surrogate models for solutions
and solve nonlinear optimization problems even for linear problems.
While neural networks are universal approximators, 
a collection of neural networks does not form a linear space.
Thus, the error estimates and convergence of deep learning methods will significantly differ from those of the LS-FEM.
The analysis of the least-squares with a linear space (e.g., LS-FEM in \citep{BocGun16,BraSch70}) cannot be applied directly.

There are two sources of errors for the deep learning algorithms for PDEs: mathematical formulation and optimization/training.
The present work focuses on analyzing the errors from mathematical formulations.
% of deep learning methods. 
The analysis of optimization/training errors is deferred to future work. 
Specifically, we study the problem of the error estimates of neural network solutions that minimize the least-squares type loss functionals.  
This problem has been investigated in 
\citep{Mishra-pinn-2023,Mishra-pinn-inv-2022} for both linear and nonlinear equations from fluid dynamics, in \citep{ShinDK20} for a discrete loss functional for linear elliptic and parabolic equations, and in \citep{SirSpi18} for a continuous loss functional 
for quasilinear parabolic PDEs. 
More related works are summarized in Section \ref{ssec:limit-rem}.

We consider a general framework regardless of the type of equations for linear problems, which include hyperbolic, elliptic, and parabolic equations.  Moreover, we consider two types of loss functionals.
%, we analyze two least-squares formulations of the loss functionals.
Type 1 is based on residuals of the strong form, and Type 2 is based on residuals of the weak form. 
% using neural networks.  
The two types of loss functionals are related to some existing methods. 
When the feasible solutions are in  
the space of finite elements, the discrete Type 1 
is known as 
least-squares collocation methods \citep{BocGun98}.
When the feasible solutions are neural networks (which we consider), Type 1 is known as physics-informed neural networks (PINNs) \citep{raissi2017machine}
and Type 2 is known under the names variational PINNs \citep{kharazmi2019variational}, VarNets \citep{khodayi2019varnet} and  \textit{hp}-variational PINNs (\textit{hp}-VPINNs)
\citep{kharazmi2020hpvpinns}.
Type 1 loss requires smooth activation functions in the networks,
while Type 2 loss does not 
as its formulation is based on variational forms of PDEs (variational residuals)
allowing non-smooth networks.

%\subsection{Contributions}
Our main contributions are summarized as follows.
%\begin{itemize}
%\item
\textit{First}, 
we establish an abstract framework 
in analyzing the discrete loss functional
whose formulation is based on 
the residuals of the linear problem \eqref{eq:linear-generic-bvp}
in strong (type 1) and weak (type 2) forms.    
% 
%
%\item 
\textit{Second}, we derive a priori and a posterior error estimates for 
both continuous and discrete loss formulations of both types: residuals in strong forms  (Sections \ref{sec:pinn-error} and \ref{sec:effect-discretization-loss}) and those in weak forms (Section \ref{sec:hp-vpinn}).
%
%\item 
\textit{Third}, under suitable assumptions, 
we establish strong convergence in the underlying predefined topology.
See Theorem~\ref{thm:convg-continuousPINNs} (continuous type 1 loss functional), Theorem~\ref{thm:convg-discretePINNs} (discrete type 1 loss functional),
and Theorem~\ref{thm:convg-hpVPINNs} (continuous type 2 loss functional).
% \item With the Rademacher complexity, we explain why the residual minimization work in high dimensions. The Rademacher complexity often grows logarithmly with the dimension (see e.g. \citep{ma2019barron}). 
%
%\item 
We also present  three examples and  validate the assumptions in our abstract framework
in \ref{sec:elliptic-advection-fractional}, including elliptic, advection-reaction, and fractional diffusion equations.  

%\item  
Our framework can be applied to the recently developed discrete/continuous PINNs and discrete/continuous \textit{hp}-VPINNs.
%\end{itemize}
Moreover, the proposed framework is more general than those in \citep{Luo2020,ShinDK20} 
where properties of considered equations are used. %
We note that in \citep{Mishra-pinn-2023,Mishra-pinn-inv-2022}, 
the stability of solutions is used 
while no convergence is guaranteed 
when both the number of parameters in the neural networks and sampling points go to infinity. 

\subsection{Assumptions}
\begin{figure}[!ht]
	\centering
	\includegraphics[width=10cm]{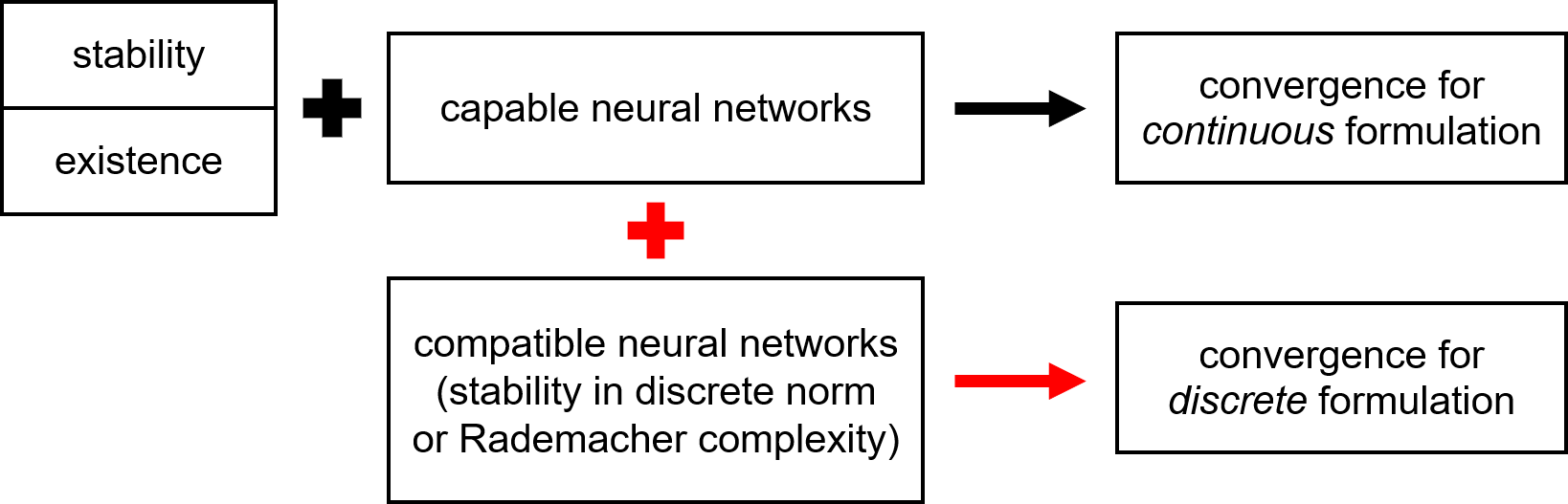}
	\caption{The sketch of critical concepts for convergence of continuous and discrete residual minimization using neural networks. Here we assume a zero training error. 
	}
	\label{fig:my_label}
\end{figure}

We have made the following three fundamental assumptions:  
(a) (\textit{stability}) relations among graph norms defined by the linear operators and the Sobolev or H\"older  norms of solutions (as in least-squares finite element methods, e.g., Assumption \ref{assum:norm-equivalence-generic} and \cite{BocGun16}), 
(b) (\textit{existence and uniqueness})
existence of a convergent sequence to solutions in such norms (Assumption \ref{assum:solution-existence}),
and (c) (\textit{compatible networks for discrete formulation}) uniformity of the discretization/projection
errors of continuous norms of residuals for minimizers, where the uniformity lies either in the stability concerning discrete norms (Assumption \ref{assum:discrete-norm-lower-bound}) or in the Rademacher complexity  (Assumption \ref{assum:Rademacher}).  The first two assumptions are somewhat standard as in most numerical methods and lead to the convergence of continuous formulations in Theorem \ref{thm:convergence-nn-linear-prior-posteriori-generic}. The third is the most essential for the convergence of discrete formulations. The violation of such an assumption may lead to no convergence, i.e., no accuracy may be gained away from the training points; see Example \ref{exm:counter-exm-poisson-no-converg}.   
In Figure~\ref{fig:my_label}, we summarize
the relations between the fundamental assumptions and convergence.

\subsection{Limitations and remarks}
\label{ssec:limit-rem}
% There have been wide applications of  least-squares formulations  using neural networks for partial differential equations,  with  great empirical success; see, e.g.,  
%\citep{BerNys18,EYu18,khoo2019solving,Lagaris_98_ANN-ODE-PDE,Lagaris_00_ANN-Irregular,LiTWL19,liaoMing19,mao2020physics,NabMei19,pang2019fpinns,raissi2018hidden,raissi2017machine,RaiPK19,SamANe20,SirSpi18,zhang2019quantifying}.
Our work is the first step toward understanding and developing neural network algorithms for forward and inverse problems.
There are aspects and issues of these algorithms that deserve in-depth investigations, e.g., choices of energy-based or residual-based loss functionals, choices of penalty parameters for optimizing various objectives,
and efficient methods for multiscale problems and nonlinear problems, especially in high dimensions. 
For example, the choice of optimum penalty parameters to balance all the residuals from different components to achieve better accuracy and convergence was discussed in \citep{BraSch70} using high-order finite element methods for elliptic problems. 
Therein some penalty parameters have been introduced to balance the two terms  $\norm{\Delta v}_{L^2(\Omega)}$   and $\norm{v}_{L^2(\partial \Omega)}$ (Dirichlet boundary) for all $v$ in a finite element space to derive optimal convergence.   
Note that the parameter can be determined via a Bernstein-type inequality.  However, such inequalities are not yet available for most neural networks, especially for deep neural networks. Some empirical study has been made, e.g., in  \citep{wang2020whenwhy}. 
% Thus, these topics are not discussed as they are beyond the scope of the present work. 
For multiscale problems, the constant in the lower bound of norm relations can be small and may lead to little use of the theory. %
Also, many intrinsic physical properties are not explicitly included in the loss functionals, such as conservation laws. 
These aforementioned topics are beyond the scope of the paper and should be considered case-by-case. 
In addition, our framework does not accommodate nonlinear problems. Some attempts have been made in \citep{Mishra-pinn-inv-2022} using the stability of solutions while more efforts are required with less restrictive assumptions; see also \citep{SirSpi18}.   
Lastly, we do not discuss the mixed formulations, which transform high-order PDEs into systems of first-order equations; see e.g., in  \citep{CaiCLL20}.

We remark that after the first version of this work \citep{shin2020error} was posted on arXiv, several related and follow-up works employed the presented abstract framework and further investigated error analysis of residual minimization. 
For example,
in \citep{gazoulis2023stability}, the concept of gamma convergence from the variation of calculus is adopted to further reduce the theoretical requirement in the smoothness of the exact solution. 
In 
\citep{Mishra_dispersive:2021}, stability is also applied 
to nonlinear equations, extending the framework in 
\citep{Mishra-pinn-2023,Mishra-pinn-inv-2022}.
The work of \citep{tang2023pinns} utilizes the error estimate for continuous residual minimization (Theorem~\ref{thm:convergence-nn-linear-prior-posteriori-generic}) to develop an adaptive sampling strategy for PINNs. 
Regarding the compatibility of neural networks,
\citep{JiaoLaiLie22} calculated the Rademacher complexity of PINNs
by focusing on the two-layer $\text{ReLU}^3$ neural networks,
which complements the presented discrete residual minimization formulation in Section~\ref{sec:prototype-PINNs} and provides a rate of convergence.
In \citep{doumeche2023convergence}, the compatibility was imposed via Sobolev-type regularization, similar to the one used in \citep{ShinDK20}.
Also, the stability was used to provide a convergence result for PINNs.
Similarly, in \citep{Wu23}, the compatibility was used by means of a Lipschitz regularization and provided a convergence result via stability for elliptic interface problems. 

Extensions to ill-posed problems with conditional stability (see e.g. \cite{Kabanikhin2008,DahMS23-leastsquare-pde} for definition) are possible. Conditional stability is an extension of stability in the sense that adding a regularization term will stabilize the problem if it is not well-posed/stable. 
For linear problems such as Poisson equations, heat equations, and wave equations, condition stability has been well established for data assimilation problems, see 
\cite{burman2016stabilised,DahMS23-leastsquare-pde} for examples.  We briefly discuss the conditional stability in Remark \ref{rem:conditional-stability} while we will not discuss the application of conditional stability in this work.

There is another NN-based approach using energy minimization.
The Deep Ritz Method, \cite{EYu18}, is a popular such approach.
While this is not the focus of the present work, for completeness, we briefly review some relevant works. 
A prior generalization analysis was given in \citep{lu2021priori,lu2022priori}
where the Barron spectral space was employed for the compatibility.
\cite{pmlr-v190-muller22a} derived an error estimate where the Friedrich inequality was used for the compatibility
and C\'ea's lemma provided the stability.
A rate of convergence was studied in \citep{Duan22} where the compatibility was handled by the Rademacher complexity
and the stability was by C\'ea's lemma.
\cite{dondl2022uniform} provided  abstract gamma convergence.

%% file: setup.tex
%%%%%%%%%%%%%%%%%%%%%%%%%%%%%%%%%%%%%%%%%%%%%%%%%5
\section{Mathematical Setup and Preliminaries} \label{sec:setup}
%%%%%%%%%%%%%%%%%%%%%%%%%%%%%%%%%%%%%%%%%%%%%%%%%5

Let $A:\, X\to Y$ and 
$B: X \to Z$ be linear operators,
where 
$(X,\|\cdot\|_X)$, $(Y,\|\cdot\|_Y)$, $(Z,\|\cdot\|_Z)$ are Banach spaces.
Let $\Omega, \Gamma$ be two subsets of $ \Real^d$
and consider the following linear problem %\eqref{eq:linear-generic-bvp}.
\begin{equation}\label{eq:linear-generic-bvp} 
    \begin{split}
        A[u](x) \hspace{0.2cm} &= \hspace{0.2cm}f(x),\quad x\in \Omega, \\
        B[u](x) \hspace{0.2cm} &= \hspace{0.2cm} g(x), \quad x \in \Gamma,
    \end{split}
\end{equation}
where $f \in Y$ and $g \in Z$. 
For simplicity, we  write 
$A[u]$ as $Au$ and similarly for $B$.  
We will impose further assumptions on the operators $A$ and $B$ shortly.

If the problem represents a linear PDE, we usually consider that 
$\Omega$ is a bounded domain (open and connected) with Lipschitz continuous boundary  
$\Gamma= \partial\Omega$.   The domain $\Omega$ could be unbounded as in the exterior problems, or the domain $\Gamma$ could be unbounded (see \ref{sec:ssec-fractional-laplacian}). For computational purposes, it is often convenient to assume bounded computational domains.  
The choices of the spaces depend on the considered problems, see \ref{sec:elliptic-advection-fractional} for several examples. 

\begin{rem}
For   simplicity,
we consider the case where the number of operators is two, 
and $u$ is a scalar real-valued function. 
% two operators $A$ and $B$ 
% and the scalar real-valued function $u$. 
However, the framework developed in the following sections can accommodate multiple operators and vector-valued functions (systems of equations), e.g., first-order differential linear systems, which can also represent elliptic problems, e.g. in \cite{BocGun98}. 
\end{rem}

\subsection{Solution of the problem and norm relations}
In this section, we present two key assumptions for our analysis.
The first assumption  is on the existence of a solution to the problem \eqref{eq:linear-generic-bvp}; see Assumption \ref{assum:solution-existence}. The second assumption is on the relations of graph norms associated with the linear problem \eqref{eq:linear-generic-bvp}; see Assumption \ref{assum:norm-equivalence-generic}.

We first define our notion of a solution to
\eqref{eq:linear-generic-bvp}, which is the underlying solution we want to approximate. 
\begin{defn} \label{def:solution}
 Let $(V,\|\cdot\|_V)$ be a Banach space and $(X,\|\cdot\|_V)$ is a dense subspace of $(V,\|\cdot\|_V)$.  
    An element $u^* \in V$
    is said to be a solution to 
    \eqref{eq:linear-generic-bvp}
    if there exists 
    a sequence $\{u_k^*\}$ in $X$
    such that 
    \begin{align*}
        \lim_{k\to \infty} \|u_k^* - u^*\|_V &= 0, \qquad 
        \lim_{k \to \infty} 
        \|Au_{k}^* - f\|_Y
        + \|Bu_{k}^* - g\|_Z = 0.
    \end{align*}
\end{defn}

We note that since a solution is not necessarily in $X$, strictly speaking,
$Au^*$ is not well-defined.
However, under suitable assumptions,
the linear operators can be uniquely extended to $V$. 
The following theorem provides a condition under which the extension exists.
Hence, $Au^*$ could be well-defined.

\begin{thm}[Bounded Linear Transformation Theorem \citep{reed2012methods}] \label{thm:BLT-theorem}
    Suppose $(X,\|\cdot\|_V)$ is a dense subspace of $(V,\|\cdot\|_V)$
    and let $A$ be a bounded linear operator from   $(X,\|\cdot\|_V)$ to $(Y,\|\cdot\|_Y)$.
    Then,
    there exists a unique extension $\tilde{A}$ from $(V,\|\cdot\|_V)$ to $(Y,\|\cdot\|_Y)$
    of $A$. 
    That is,
    $Av = \tilde{A}v$ 
    for all $v \in X$
    and $\|A\| = \|\tilde{A}\|$,
    where $\|\cdot\|$ is the operator norm.
\end{thm}

Throughout the paper, the operators are understood up to extensions if needed.

\begin{assu}[\textbf{Existence}] \label{assum:solution-existence}
    There exists a solution $u^* \in V$
    to \eqref{eq:linear-generic-bvp}
    in the sense of Definition~\ref{def:solution}.
\end{assu}

The following  assumption plays a central role in 
our abstract framework. 
\begin{assu}[\textbf{Norm relations}]\label{assum:norm-equivalence-generic}
	Assume that the operators $A$ and $B$ satisfy 
	$\|Av\|_Y, \|Bv\|_Z < \infty$,  for all $v \in V$.
	Assume also the following norm relations: 
	\vskip -10pt 
	\begin{subequations}\label{eq:norm-equivalence-generic}
	\begin{equation}\label{eq:norm-equivalence-generic-A}
	   C_1\norm{u}_V \leq  \norm{Au}_Y +\norm{Bu}_Z, \quad \forall u \in X, 
	\end{equation}
	\vskip -15pt
	\begin{equation}\label{eq:norm-equivalence-generic-B}
    \norm{Au}_Y +\norm{Bu}_Z\leq  C_2 \norm{u}_{X}, \quad \forall  u \in X,
 \end{equation} 
 \end{subequations}
 where the positive constants  $C_1$ and $C_2$ 
	do  not depend on $u$ but on the domain and the coefficients of the operators $A$ and $B$.
As before, 	$X$ is a dense subspace of $V$.
\end{assu} 
%\vskip -20pt
This assumption on norm relations is not new and 
has been used for many least-squares formulations for numerical methods such as  finite element methods, e.g. 
\cite{BocGun16,BraSch70}. %,BraSch71,BraTho72}. 
The first norm relation (\eqref{eq:norm-equivalence-generic-A} in Assumption~\ref{assum:norm-equivalence-generic}) gives the \textit{stability/regularity} of the solution. 
The second norm relation (\eqref{eq:norm-equivalence-generic-B} in Assumption~\ref{assum:norm-equivalence-generic}) is for smooth numerical solutions   in $X$ in which we need universal approximation by neural networks, instead of in $V$; See \ref{sec:elliptic} for  Case I of an elliptic problem.
Also, the condition of \eqref{eq:norm-equivalence-generic-A} in Assumption~\ref{assum:norm-equivalence-generic}
guarantees the uniqueness of the solution in $V$
as shown in the following proposition.
\begin{prop}[\textbf{Uniqueness}]
   Let the condition of \eqref{eq:norm-equivalence-generic-A} in Assumption~\ref{assum:norm-equivalence-generic} hold.
    Then, the solution to \eqref{eq:linear-generic-bvp}
    (Definition~\ref{def:solution})
    is unique if it exists.
\end{prop}
\begin{proof}
     It suffices to consider 
     the problem with homogeneous data $f=g=0$.
     Let $u^*$ be a solution 
     to \eqref{eq:linear-generic-bvp}
     in the sense of Definition~\ref{def:solution}. 
     Then, there exists 
     a sequence  $\{u_k\}$ in $X$
     such that $\|u_k - u^*\|_V, \|Au_k\|_Y, \|Bu_k\| \to 0$.
     If
    \begin{align*}
        C_1 \norm{u_k}_V &\leq  \norm{Au_k}_Y +\norm{Bu_k}_Z 
        \to 0, \text{ then } \norm{u_k}_V \to 0.
    \end{align*}
    Therefore, $u^* = 0$,
    which shows the uniqueness. 
\end{proof}

\begin{rem}\label{rem:conditional-stability}
The stability (well-posedness) may be relaxed to conditional stability. 
If the pair $(A, B)$ is not unconditionally stable, i.e., in the sense of \eqref{eq:norm-equivalence-generic-A}, we may consider a regularization operator $L$ and the triplet 
$(A, B, L)$, where $L: X \to H$ is linear, where $H$ is a proper Banach space such that it is stable in the following sense:
there exist a positive function $\mathcal{D}:X \to \Real^{+}$
and  a non-decreasing  function
$\rho_C:\,\Real^+\cup\set{0} \to  \Real^+\cup\set{0}$ with $\lim_{t\searrow 0}\rho_C(t) = 0$ for any $C >0$, such that  for $v\in X$ with $\norm{Lv}_H\leq C$, it holds that
\[\mathcal{D}(v)\leq \rho_C(\norm{Av}_X +   \norm{Bv}_Z).\]
Here $\mathcal{D}$ is usually a norm or a semi-norm on a Banach space that $X$ can be embedded into. 
When $\mathcal{D}$ is a norm, we do not need a regularization ($L$) as the problem is well-posed.   The loss function at the continuous level can be formulated as 
\[ \norm{Av -f}_Y^2 + \norm{Bv-g}_Z^2 + \epsilon^2\norm{Lv}_H^2,\quad \epsilon \text{ is small}.\]
 The conclusion in the following sections can be modified accordingly to obtain stability and  
 convergence of the PINN formulations.
 We will limit ourselves to well-posed problems for a simple presentation. 
We refer to Section 3 of \cite{DahMS23-leastsquare-pde} and the references therein for examples. 
\end{rem}
 
\subsection{Loss Functionals of Residual Minimization}
We present four  loss functionals for 
the residual minimization (RM).

1) (Discrete RM) Given $\{x_i^r, f(x_i^r) \}_{i=1}^{M_r}$, $\{x_i^b, g(x_i^b)\}_{i=1}^{M_b}$
where $\{x_i^r\} \subset \Omega$, $\{x_i^b\} \subset \Gamma$,
    the discrete loss functional for RM is defined by
    \begin{equation*}
        \mathcal{J}^M(u) = \frac{1}{M_r}\sum_{i=1}^{M_r} (f(x_i^r) - A[u](x_i^r))^2 + \frac{1}{M_b}\sum_{i=1}^{M_b} (g(x_i^b) - B[u](x_i^b))^2,
    \end{equation*}
    where $M = (M_r, M_b)$.  See e.g., \cite{raissi2017machine}. 
  % This formulation is used in \cite{raissi2017machine}. 
    
2)  (Continuous RM)
    The continuous loss functional for RM is defined by
  \begin{equation*}
      \mathcal{J}(u) = \|f - Au\|_{L^2(\Omega)}^2 + \|g - Bu\|_{L^2(\partial \Omega)}^2.
  \end{equation*} 

Let $\Omega$ be a bounded domain in $\mathbb{R}^d$
and $\{\Omega_k\}_{k=1}^K$ be a partition of $\Omega$.
For each $k$, let $\{\Phi_{k,s}\}_{s \ge 1}$ be
a  complete orthonormal basis in $L^2(\Omega_k)$ or $H^1_0(\Omega_k)$.
Then, by defining $\Phi_{k,s}(x) = 0$ for all $x \notin   \Omega_k$ and  $k,\,s$, 
 $\{\Phi_{k,s}\}$ is a complete orthonormal basis of $L^2(\Omega)$ or $H^1_0(\Omega)$.

3) (Discrete hp-VRM)  Let $\{x_{i}^{k}, f(x_i^{k}) \}_{i=1}^{M_{k}}$ and $\{x_i^b, g(x_i^b)\}_{i=1}^{M_b}$ be training points
    and 
    $\{\Phi_{k,i}\}_{i=1}^{N_k}$ be a set of test functions in $L^2(\Omega_k)$ for $1\le k \le K$.
    Let $\bm{N} = (N_1,\cdots, N_K)$ and $\bm{M} = (M_1,\cdots, M_K, M_b)$.
    Then, a version of the discrete loss functional for hp-variational RM (hp-VRM) is given by
    \begin{eqnarray*} 
            \mathcal{J}^{\bm{M},\bm{N}}(u)  = \sum_{k=1}^K \sum_{s=1}^{N_k} \left(\sum_{i=1}^{M_{k}} w^{k}_i (f(x_{i}^{k}) - A[u](x_{i}^{k}))\Phi_{k,s}(x_{i}^{k}) \right)^2 + \frac{1}{M_b}\sum_{i=1}^{M_b} (g(x_i^b) - B[u](x_i^b))^2. 
    \end{eqnarray*}
    By applying integration by parts (repeatedly if needed), 
    a different but equivalent formulation of RM can be obtained where the regularity requirement for the numerical solution and the activation functions is weakened; see e.g.,  \cite{kharazmi2020hpvpinns,kharazmi2019variational}.

4) (Continuous hp-VRM)
    The continuous loss functional for hp-VRM is given by
   \begin{equation*}
        \mathcal{J}^{\bm{N}}(u) = \sum_{k=1}^K \sum_{s=1}^{N_k} |\langle f-Au, \Phi_{k,s} \rangle_{L^2(\Omega_k)}|^2 +\|g - Bu\|_{L^2(\partial \Omega)}^2.
   \end{equation*}

%%================================================
%================================
\subsection{Feed-forward Neural Networks}
Throughout the paper,  any neural networks can be applied in our framework 
as long as they satisfy the universal approximation theorem  required in the space of $X$.  For presentation, we assume that the feasible solutions are 
feed-forward neural networks unless explicitly stated  otherwise.

Let us first review some of known universal approximation 
theorems for feed-forward neural networks.
For a positive integer $L$, 
a $L$-layer feed-forward neural network is a function $\mathcal{R}[\theta]:\mathbb{R}^{n_0}\mapsto \mathbb{R}^{n_L}$  
that is recursively defined by
\begin{equation*}
    \mathcal{R}[\theta](x) = z^L(x), \quad z^l(x) = W_l\phi(z^{l-1}(x)) + b_l, \quad 
    2 \le l \le L, \quad 
    z^1(x) = W_1x + b_1.
\end{equation*} 
Here \(W_l\in\mathbb R^{n_{l}\times n_{l-1}}\) is the $l$-th layer weight matrix  and \(b_l\in\mathbb R^{n_l}\) is the $l$-th layer bias vector. The activation 
$\phi(x)$ is applied elementwise.   
The collection of network parameters is $\theta = \set{(W_1, b_1), \dots, (W_L, b_L)}$.
The architecture of the network is represented by the vector
$\vec{\bm{n}} = (n_0,\dots, n_L)$.
With the fixed architecture $\vec{\bm{n}}$,
the set of all possible network parameters 
is denoted by \footnote{It is possible to use bounded weights as the resulting network is still a universal approximator. We will not consider this for simplicity.}
\begin{equation} \label{def:NN-parameter-space}
\Theta(\vec{\bm{n}}) = \left\{ \{(W_j, b_j)\}_{j=1}^L : W_j \in \mathbb{R}^{n_j \times n_{j-1}}, b_j \in \mathbb{R}^{n_j},\,j=1,2,\ldots,L. \right\}. 
\end{equation}
%
%%%%%%%%%%%%%%%%%%
%
%
Let $\{\vec{\bm{n}}_n\}_{n \ge 1}$
be a sequence of network architectures
such that $\vec{\bm{n}}_n \le \vec{\bm{n}}_{n+1}$ for all $n$,
where the vector inequality is understood entry-wise.
We then define its corresponding sequence of neural network classes by
\begin{equation} \label{def:NNs}
    \mathfrak{N}_{\theta,n} = \left\{ \mathcal{R}[\theta](x) : \theta \in \cup_{\vec{\bm{v}} \le \vec{\bm{n}}_n}\Theta(\vec{\bm{v}}) \right\}.
\end{equation}
% {\color{red}where neural networks of interest are the ones having finite $\|\cdot\|_X$-norm.}
% This implicitly restricts weights and biases to be bounded.
%
By the construction, we have $\mathfrak{N}_{\theta,n} \subset \mathfrak{N}_{\theta,n+1}$.
An element of $\mathfrak{N}_{\theta,n}$
is simply denoted by $u_{\bN, n}$,
where $\bN$ stands for neural network (NN).

Next, we make the following assumption
on the sequence of network classes, 
which guarantees
universal approximation.
\begin{assu}[\textbf{Uniform NN approximation of elements in $X$}]\label{assum:approximation-abilty-X}
%    \textcolor{blue}{Let $(X,\|\cdot\|_V)$ be a dense subspace of $(V,\|\cdot\|_V)$.}
%    Let $X$ be a dense subspace of $V$.
    % and $V_0$ be the set defined 
    % in \eqref{def:V0}.
    There exists a sequence of neural network classes
    such that $\mathfrak{N}_{\theta,n} \subset \mathfrak{N}_{\theta,n+1}$
    and   \,
    {$X \subset \overline{\bigcup_n \mathfrak{N}_{\theta,n}}$} in the toplogy of ($X,\|\cdot\|_X)$.  
%For any $v\in X$, there exists a sequence  in  $\overline{\bigcup_n \mathfrak{N}_{\theta,n}}$ such that 

\end{assu}

In literature, Assumption~\ref{assum:approximation-abilty-X} is proved for various spaces $X$. For example,  the work of \cite{Guhring_19}
shows that 
for any 
$\epsilon \in (0, 0.5)$,
$p \in [1,\infty]$,
and $s \in [0, 1]$, 
there exists a deep rectified linear unit (ReLU) network architecture $\vec{\bm{n}}$
such that for any $f \in W^{n,p}((0,1)^d)$ with $\|f\|_{W^{n,p}}$ being bounded for all $f$ and $n \ge 2$,
there exists a neural network with parameters $\theta \in \Theta(\vec{\bm{n}})$ with 
%\begin{equation*}
 $   \|\mathcal{R}[\theta] - f \|_{W^{s,p}} \le \epsilon$.
%\end{equation*}
Here the network architecture
depends only on $d, n, p, s, \epsilon$ and the uniform bound of $f$. 
Also, the work of \cite{mhaskar1997neural}
showed a similar result
for two-layer networks with smooth activation functions.

%%%%%%%%%%%%%%%%%%%%%%%%%%%%%%%%%%%
% end of setup, section 
%%%%%%%%%%%%%%%%%%%%%%%%%%%%%%%%%%%

%% file: crm.tex
\section{Continuous Residual Minimization}\label{sec:pinn-error}
A goal  of the residual minimization is to approximate the solution to Equation  \eqref{eq:linear-generic-bvp}  
by solving the optimization problem 
\begin{equation}\label{eq:pinn-opt-problem}
\inf_{v\in \mathfrak{N}_{\theta,n}\cap X} \mathcal{J}_\tau(v),
\end{equation}
where the loss functional $\mathcal{J}_\tau(v)$ is defined by \footnote{In general, it can be 
$\mathcal{J}_\tau(v ) = \tau_r\norm{f-Av}_Y^{p} + \tau 
\norm{g-Bv}_Z^{p}$, where $\tau_r,\tau>0$. In this work, we focus on $\tau_r=1$.}
\begin{equation}\label{eq:loss-functional}
\mathcal{J}_\tau(v ) = \norm{f-Av}_Y^{p} + \tau 
\norm{g-Bv}_Z^{p}.
\end{equation}
Here $\tau > 0$, $p \ge 1$, and $\tau$ is a fixed parameter that weighs
the discrepancy of the boundary 
in the loss. Also, we assume that $f\in Y$ and $g\in Z$ and otherwise a mollifier may be applied; see e.g., in Remark \ref{rem:non-smooth-data}. 

\begin{rem}
In practice, 
it has been empirically shown that 
the choice of $\tau$ 
significantly influences 
the training of neural networks \cite{Lagaris_00_ANN-Irregular,wang2020understanding,wang2020whenwhy}. It is observed that $\tau$ depends on the network but may be only moderately large, in which   $\tau$ will not  affect 
 our analysis up to some constant. % as long as $\tau \ge 1$.
 \end{rem}

\begin{rem}
    For simplicity, 
    we often assume 
    the existence of the global minimizer 
    of \eqref{eq:pinn-opt-problem}. For example, 
    the global minimizer exists when a two-layer network is used for solving PDEs; see \cite{Luo2020} for details.
%    throughout the paper.
    This assumption can be relaxed 
    by introducing 
    appropriate 
    quasi-minimizers.
\end{rem}

The choices of $Y$ and $Z$ are important to design
the loss functionals; see \ref{sec:elliptic-advection-fractional} for examples. A guideline for choosing 
$Y$ and $Z$ is Assumption \ref{assum:norm-equivalence-generic}, which depends on the considered problems and the metric of accuracy (the $V$-norm).  
For PDEs,  we may have $\Gamma=\partial\Omega$ and some typical choices for $Y$ and $Z$ are 
  as follows:
\begin{itemize}
    \item (Banach spaces) 
    $Y=L^p_\rho(\Omega)$ and $Z=L^p_{\rho_b}(\partial\Omega)$, where $p\geq 1$, $\rho: \Omega\to [0,+\infty)$, $\rho \in L^1(\Omega)$, $\rho_b : \partial\Omega \to [0,+\infty)$, and $\rho_b\in L^1(\partial\Omega)$.
    \item (Banach spaces) 
    $Y=L^\infty(\Omega)$ and $Z=L^\infty(\partial\Omega)$.
%    ,\quad $p=2$.
   %  \item (Banach spaces)      $Y=C(\Omega)$ and $Z=C(\partial\Omega)$.
    \item (Hilbert spaces) $Y= H^{k}(\Omega)$, $k=-2,-1,0, 1,2,\cdots$, and  $Z=L^2_{\rho_b}(\partial\Omega)$.
\end{itemize}

\subsection{Error estimates}
%If the hypothesis class were the entire space $X$, 
%Assumption~\ref{assum:solution-existence}
%directly results in the zero loss.

The convergence of (quasi-)minimizers of the loss functional depends on the  
neural network classes $\{\mathfrak{N}_{\theta,n}\}_{n\ge 1}$. 
%converges or not depends on  $\{\mathfrak{N}_{\theta,n}\}_{n\ge 1}$.
The following proposition shows that 
the universal approximation property 
(Assumption~\ref{assum:approximation-abilty-X})
{and 
 the norm relation of \eqref{eq:norm-equivalence-generic-B} 
 in  Assumption~\ref{assum:norm-equivalence-generic}}
are sufficient for the convergence of the loss. % a prior and a posterior error estimates.
% Under some mild conditions 
% on hypothesis classes, 
% can also achieve the zero loss, 
% which is shown in the following proposition.
\begin{prop} \label{prop:loss-to-zero}
    Suppose 
    Assumptions~\ref{assum:solution-existence},~\ref{assum:approximation-abilty-X} and
    {
    \eqref{eq:norm-equivalence-generic-B}
    of Assumption~\ref{assum:norm-equivalence-generic} hold.}
    % $\mathfrak{N}_{\theta,n} \subset \mathfrak{N}_{\theta,n+1}$ and 
    % $\overline{\bigcup_n \mathfrak{N}_{\theta,n}} \cap X = X$.
    For a fixed $\tau > 0$, 
    let $u_{\bN,n}^\tau$ be a quasi-minimizer \footnote{It means that  $\mathcal{J}_\tau(u_{\bN,n}^\tau) \leq \inf_{v \in \mathfrak{N}_{\theta,n}\cap X} \mathcal{J}_\tau(v) + \epsilon$, where $\epsilon\geq 0$ is small. }
  % \inf_{v\in }  \mathcal{J}_\tau(v) +$}
     of the loss \eqref{eq:loss-functional}.
    Then,
   % \begin{equation*}
    {$\lim_{n\to \infty} \mathcal{J}_\tau(u_{\bN,n}^\tau) = 0$.}
    %\end{equation*}
\end{prop}
{
\begin{proof}
	Let $v^*$ be a solution to \eqref{eq:linear-generic-bvp}
	and $\{v_k^*\}$ be its corresponding 
	sequence in $X$ (Definition~\ref{def:solution}).
	Let $\{\epsilon_k\}$ be a positive decreasing sequence that converges to 0.
	For $k \gg 1$,
	it follows from Assumption~\ref{assum:approximation-abilty-X}
	that there exists an integer $n_k$ 
	and a neural network 
 $u_{n_k} \in \mathfrak{N}_{\theta,n_k}$
 % $u_{n_k} \in \mathfrak{N}_{\theta,n_k} \cap X$
	such that $\|v_k^* - u_{n_k}\|_X \le \epsilon_k$\footnote{This implies $u_{n_k} \in \mathfrak{N}_{\theta,n_k}\cap X$ as $\|u_{n_k}\|_X \le \epsilon_k + \|v_k^*\|_X < \infty$.}.
	Observe that Equation~\eqref{eq:norm-equivalence-generic-B} gives
	\begin{align*}
		\mathcal{J}_\tau(u_{n_k})
		&\le (\|f-Av_k^*\|_Y + C_2\epsilon_k)^p + 
		\tau  (\|g - Bv_k^*\|_Z + C_2 \epsilon_k)^p.
	\end{align*}
	Let $\{u_{\bN,n}^\tau\}$ be a sequence of quasi-minimizers of the loss functional in $\mathfrak{N}_{\theta,n}\cap X$,
	i.e., $\mathcal{J}_\tau(u_{\bN,n}^\tau)
	\le \inf_{v \in \mathfrak{N}_{\theta,n}\cap X} \mathcal{J}_\tau(v) + \delta_n$,
	where $\delta_n \to 0$ as $n \to \infty$.
	Observe that $\mathcal{J}_\tau(u_{\bN,n_k}^\tau) \le \mathcal{J}_\tau(u_{n_k}) + \delta_{n_k}$.
	By letting $k \to \infty$,
	we have 
	$\lim_{k\to \infty} \mathcal{J}_\tau(u_{\bN,n_k}^\tau)  = 0$. 
	For $m_1 < m_2$,
	since 
	$\mathfrak{N}_{\theta,m_1}\subset\mathfrak{N}_{\theta,m_2}$,
	we have $\mathcal{J}_\tau(u_{\bN,m_2}^\tau)
	\le \mathcal{J}_\tau(u_{\bN,m_1}^\tau) + \delta_{m_2}$.
	Hence, 
	we can conclude that 
	$\lim_{n\to \infty} \mathcal{J}_\tau(u_{\bN,n}^\tau) = 0$\footnote{For any $\epsilon > 0$, there exists $K$ such that $\mathcal{J}_\tau(u_{\mathbb{N},n_k}^\tau) \le \epsilon/2$
		for all $k \ge K$.
		Also, there exists $N$ such that 
		$\delta_n \le \epsilon/2$ for all $n \ge N$.
		By choosing $\hat{N} = \max\{n_K, N\}$,
		we have 
		$\mathcal{J}_\tau(u_{\mathbb{N},n}^\tau) \le \epsilon$
		for all $n \ge \hat{N}$.
	}.
\end{proof}
}

As shown in Proposition~\ref{prop:loss-to-zero}, a zero-loss can be achieved 
if the neural network classes $\mathfrak{N}_{\theta,n}$
can capture an approximating sequence in $X$ (in the sense of Definition~\ref{def:solution}),
{and the norm relation of \eqref{eq:norm-equivalence-generic-B} holds.}
However, the convergence of the loss 
does not necessarily imply 
the convergence of the (quasi-)minimizers to the solution to the governing equation \cite{ShinDK20}.
 
Next, we present  
\emph{a priori} and \emph{a posterior} error estimates for  optimization of \eqref{eq:pinn-opt-problem}. 

%%%%%%%%%%%%%%%%%%%%%%%%%%%%%%%%%%%%%%%%%%%%%%%%%

  \begin{thm}[\textbf{Error estimates for continuous RM}] \label{thm:convergence-nn-linear-prior-posteriori-generic} 
  	Let Assumptions
\ref{assum:solution-existence}
and 
\ref{assum:norm-equivalence-generic} hold. 
Let $u_{\mathbb{N},n}^{\tau} \in \mathfrak{N}_{\theta,n} \cap X$ be a solution to the minimization problem  \eqref{eq:pinn-opt-problem}
and $u^*$ be the solution to \eqref{eq:linear-generic-bvp}
from Assumption~\ref{assum:solution-existence}.
Then for $\tau\geq 1$,  the following \emph{a posterior estimation} holds:
\begin{equation}\label{eq:a-posteriori-error-generic}
\norm{u_{\mathbb{N},n}^{\tau} - u^*}_{V}\leq  C_1^{-1}2^{\frac{p-1}{p}} \left(\mathcal{J}_{\tau}(u_{\mathbb{N},n}^{\tau})\right)^{\frac{1}{p}}.
\end{equation} 
(\emph{A priori estimate}) Also, for any $\epsilon > 0$, there exists 
$u_{\epsilon}^* \in X$ such that 
\begin{equation}\label{eq:a-priori-error-generic}
        \|u_{\mathbb{N},n}^{\tau} - u^*\|_V   \leq    2^{\frac{p-1}{p}}(1+\tau)^{1/p}	C_1^{-1} \big ( C_2\inf_{w\in \mathfrak{N}_{\theta,n}\cap X}\norm{w -u_\varepsilon^*}_X  + \epsilon
        \big).
\end{equation}
Here the constants $C_1$ and $C_2$ are  defined in the norm relations 
\eqref{eq:norm-equivalence-generic}.
\end{thm}
\begin{proof}
	Since $\tau \ge 1$, it then follows from 
	%this limit     and 
	the condition \eqref{eq:norm-equivalence-generic-A}
    that 
    \begin{equation}\label{eq:through-norm-relation-lower-bound}
    C_1\norm{u^*- w }_V \leq      
    2^{\frac{p-1}{p}}  \big(\mathcal{J}_{1}(w\big) )^{\frac{1}{p}}\leq     
     2^{\frac{p-1}{p}}   \big(\mathcal{J}_{\tau}(w\big) )^{\frac{1}{p}}, \qquad \forall w\in X \subseteq V.
\end{equation} 
Letting  $w=u_{\mathbb{N},n}^{\tau}$ leads to the estimate \eqref{eq:a-posteriori-error-generic}.
By \eqref{eq:through-norm-relation-lower-bound}, we have 
\begin{eqnarray}\label{eq:bound-optimization}
	C_1   \|u_{\mathbb{N},n}^{\tau} - u^*\|_V & \leq&  2^{\frac{p-1}{p}}  \big(\mathcal{J}_{\tau}(u_{\mathbb{N},n}^{\tau}\big) )^{\frac{1}{p}}
	=  2^{\frac{p-1}{p}}\inf_{w\in \mathfrak{N}_{\theta,n}\cap X}    \big(\mathcal{J}_{\tau}(w\big) )^{\frac{1}{p}} \notag 
 \notag \\
&\leq&   2^{\frac{p-1}{p}}(1+\tau)^{1/p}\inf_{w\in \mathfrak{N}_{\theta,n}\cap X}  (\norm{Aw - f}_Y 
	+\norm{Bw -g}_Z).
\end{eqnarray}
For $\epsilon > 0$,
let $u_{\epsilon}^* \in X$ such that
$\|Au_{\epsilon}^* - f\|_Y + \|Bu_{\epsilon}^*-g\|_Z < \epsilon$
% Let $ u_\varepsilon^*$ be a mollified $u^* \in X$
% satisfying 
% $\norm{ Au_\varepsilon^* -A u^*}_Y +  \norm{B u_\varepsilon -u^*}_Z < \epsilon$
from Assumption~\ref{assum:solution-existence}.
By the triangle inequality, we have 
 \begin{eqnarray} \label{eq:tri-error-mollif-appr}
 &&\norm{A w -f}_Y +  \norm{B w -g}_Z \notag \\
 &\leq&  	(\norm{Aw -Au_\varepsilon^*}_Y +  \norm{Bw  -Bu_\varepsilon^*)}_Z)  +
 \norm{ Au_\varepsilon^* -f}_Y +  \norm{B u_\varepsilon^*  -g }_Z %\notag \\ &\le & 
 \notag \\
 &<&  C_2\norm{w - u_\varepsilon^*}_X
 + \epsilon,
%  \mathcal{J}_1( u_\varepsilon^*). \qquad
 \end{eqnarray}
 where we have used the norm relation \eqref{eq:norm-equivalence-generic-B} in the last inequality.
 By combining it with \eqref{eq:bound-optimization}, we obtain   
\begin{eqnarray*}
	C_1   \|u_{\bN,n}^\tau - u^*\|_V  
%	&\leq&     2^{\frac{p-1}{p}}(1+\tau)^{1/p}\inf_{w\in \mathfrak{N}_{\theta,n}\cap X}  (\norm{Aw - f}_Y
%	+\norm{Bw -g}_Z)\\
	&\leq & 2^{\frac{p-1}{p}}(1+\tau)^{1/p} \big ( C_2\inf_{w\in \mathfrak{N}_{\theta,n}\cap X}\norm{w -u_\varepsilon^*}_X  + \epsilon
% 	+ \mathcal{J}_1( u_\varepsilon^*)
	\big),
\end{eqnarray*}
which completes the proof.
We note that if $u^*\in X$, 
one may set $u^*_{\epsilon} =u^*$ and $\epsilon = 0$.  
\end{proof}

By combining 
Theorem~\ref{thm:convergence-nn-linear-prior-posteriori-generic}
and Proposition~\ref{prop:loss-to-zero},
the convergence of the (quasi-)minimizers can be readily established.
We formally state the result below.
\begin{thm}[Convergence for continuous RM] \label{thm:convg-continuousPINNs}
    Suppose 
    Assumptions
    ~\ref{assum:solution-existence},
    ~\ref{assum:norm-equivalence-generic}, 
    and ~\ref{assum:approximation-abilty-X} hold.
    % $\mathfrak{N}_{\theta,n} \subset \mathfrak{N}_{\theta,n+1}$ and 
    % $\overline{\bigcup_n \mathfrak{N}_{\theta,n}} \cap X = X$.
    For a fixed $\tau \ge 1$, 
    let $u_{\bN,n}^\tau$ be a quasi-minimizer of the loss functional \eqref{eq:loss-functional}.
    Then,
   $ % \begin{align*}
        \lim_{n\to \infty} \|u_{\bN,n}^\tau - u^*\|_V = 0.
   $ % \end{align*}
\end{thm}

%%=============================================================

%%=============================================================
\begin{rem}
In practice, regularization is often used in the loss functional \eqref{eq:loss-functional} 
to improve the efficiency of numerical methods. 
We consider the regularization of non-smooth functionals. 
For example,  when $Y=L^1(\Omega)$ (as in  \cite{Guermond04lp}) and $Z=L^1({\partial\Omega})$,
the loss functional is not Fr\'{e}chet differentiable.
In such cases, the loss functional is hard to optimize using gradient-based methods. 
%
%In Section \ref{sec:ssec-lp-error}, we consider a regularization when $Y=L^p$, e.g. when $p=1$ as in  \cite{Guermond04lp}.
\end{rem}

\iffalse 
\begin{rem}\label{rem:approximation-error-indepedent of-u}
	The effect of discretizations of the integrals in  the functional $\mathcal{J}_\tau$ is much more complicated than the regularization considered above. The main reason is that the discretization errors of the integrals depend on the $u^*$ and $u_{\bN,n}$ while the error from regularization is independent of $u^*$ and $u_{\bN,n}$. 
\end{rem}
\fi

%%-----------------------------------------
\section{The effect of discretizing the loss functionals}
\label{sec:effect-discretization-loss}
The continuous norms (integrals) in the loss functional \eqref{eq:loss-functional} 
are discretized for simulations in finite arithmetic computers.
In this section,
we discuss the effect of 
the discretization of the loss functional \eqref{eq:loss-functional}.
%In this section, we derive error estimates of the discrete RM.
The results are based on Theorem~\ref{thm:convergence-nn-linear-prior-posteriori-generic}.
The core    is to quantify how well the discretization approximates its
corresponding continuous norm.
In order to  
characterize the relation between
the discrete norm 
and the continuous norm,
we utilize well-designed 
classes of neural networks,
%with compactness,
which guarantee
the convergence 
with respect to the number of training data samples.
The failure of using tailored classes may
lead to no convergence even if a zero training loss is achieved.
%; see  Example \ref{exm:counter-exm-poisson-no-converg}. }
%%%
In  Example \ref{exm:counter-exm-poisson-no-converg}, we show that 
even if a zero training loss is achieved, there is no convergence guarantee. 
   
\begin{exm}[Counterexample]\label{exm:counter-exm-poisson-no-converg}
	Consider the 1D Poisson equation $	-\Delta u(x) = f(x) $ on $\Omega=  (0, 1) $  with the Dirichlet boundary conditions $
	u(x) = 0 $ on $ \partial \Omega =\{0, 1\}  $. 
	Suppose that $u^*$ is the unique  classical  solution. 
	The discrete loss functional is given by
	$
	%\begin{align*}
	\mathcal{J}^M(u) = \frac{1}{M_r} \sum_{i=1}^{M_r} w_i|\Delta u(x_i) - f(x_i) |^2
	+ \frac{1}{2}\left(|u(0)|^2 + |u(1)|^2 \right),
	$ %\end{align*}
	where $M = (M_r, M_b=2)$, $w_i$'s are some quadrature weights
	and $\{x_i\}_{i=1}^{M_r}$ with $x_i = \frac{i}{M_r+1}$. % is the set of $M_r$ equidistant points on $\overline{\Omega}=[0,1]$, i.e.,  
	Let $u_{\bN}^{M}$ be a minimizer in a network class.
	Let us consider the case where $u_{\bN}^{M}(0) = u_{\bN}^{M}(1)=0$, and
	$
	\Delta u_{\bN}^{M}(x) - f(x) ={(2 \pi (M_r+1))^k}\sin(2 \pi (M_r+1) x)
	$
	{for some positive integer $k\ge 2$},
	assuming that $f$ is sufficiently smooth. 
	It then can be checked that the training loss is zero,
	while $\|\Delta u_{\bN}^{M}- f \|_{L^2(\Omega)}^2 =
	{(2 \pi (M_r+1))^{2k}}0.5$.
	Hence, 
	$\mathcal{J}(u_{\bN}^M) = {(2 \pi (M_r+1))^{2k}}0.5$
	where 
	$\mathcal{J}(u) =  \|\Delta u - f \|^2_{L^2(\Omega)} + \frac{1}{2}\left(|u(0)|^2 + |u(1)|^2 \right)$.
	This indicates that
	minimizing the discrete loss
	may not lead to 
	the convergence in the continuous norm. 
        Even no convergence in $L^2$-norm is achieved. 
	Therefore, extra conditions
	are required to guarantee convergence. 	
\end{exm}

Let us focus on the case where $Y=L^p(\Omega)$ and $Z= L^p(\Gamma )$.
The discretized loss functional of \eqref{eq:loss-functional}
is defined as  
%\vskip -10pt \noindent
\begin{equation} \label{def:discrete-pinn-loss}
	\mathcal{J}_{\tau}^{M}(v) = \sum_{i=1}^{M_r} w_i^r (f(x_i^r) - Av(x_i^r))^p + \tau\sum_{i=1}^{M_b} w_i^b(g(x_i^b) - Bv(x_i^b))^p,
\end{equation}
\vskip -10pt \noindent
where $\{(x_i^r, w_i^r)\}_{i=1}^{M_r}$ is a set of   training points and weights in $\Omega$,
$\{(x_i^b, w_i^b)\}_{i=1}^{M_b}$ is a set of   training points and weights in $\Gamma $,
and $M = (M_r, M_b)$.
When $p=2$, $w_i^r = 1/M_r$ and $w_i^b = 1/M_b$ for all $i$, 
we recover the loss function  used in \cite{raissi2017machine}.
%
%Let $Y_{M_r}$ be a subspace of $Y$ 
	Let $Y_{M_r}$ be a set of {(some special)} functions $v \in Y$ 
	such that $v$ is continuous at points $x_i^r$
	and 
	%\vskip -10pt \[    Y_{M_r} = \{ v \in Y :  \norm{v}_{Y_{M_r}} < \infty \},
	%\quad
	% \text{where}
	%\quad 
	$ \norm{v}_{Y_{M_r}} :=   \left(\sum_{i=1}^{M_r} w_i^r v^p(x_i^r)\right)^{1/p}<\infty $.  
	{
	In many cases, 
	it can be checked that $Y_{M_r}$ 
	is a restricted subset of $Y$}.
Similarly, we define $Z_{M_b}$ as {a subset} of $Z$.
% and its inner product and norm.
With these notation,
one can simply write
the discrete loss functional 
\eqref{def:discrete-pinn-loss}
as 
$\mathcal{J}^M_{\tau}(v)
= \|f-Av\|_{Y_{M_r}}^p
+ \tau \|g - Bv\|_{Z_{M_b}}^p$.

%%%%
%%%% \newpage
 
\subsection{Using discrete norm relations} \label{sec:ssec-effect-discretizing-integrals}  
Motivated by Example \ref{exm:counter-exm-poisson-no-converg}, 
we make Assumption \ref{assum:discrete-norm-lower-bound}. 
The first condition in Assumption \ref{assum:discrete-norm-lower-bound} is violated in 
Example \ref{exm:counter-exm-poisson-no-converg}, e.g., when $f=0$.

%%%%%%%%%%%%%%%%%%%%%
 
\begin{assu}\label{assum:discrete-norm-lower-bound}
	{
	For any $n \in \mathbb{N}$,
	there exist 
%    two sets $Y_c\subseteq Y$
%    and $Z_c \subseteq Z$
    two positive integers $M_r$ and $M_b$
	such that 
%    \vskip -10pt
%  
     \begin{equation} \label{eqn:discrete-norm-relations}
    \begin{split}
    \|A v  \|_{Y_{M_r}} &\ge \frac{1}{2}\|Av \|_{Y}, 
    \quad 
    \forall v \in   \set{v \in \mathfrak{N}_{\theta,n} \cap X |  Av   \in Y_{M_r}},
    \\
    \|Bv \|_{Z_{M_b}} &\ge \frac{1}{2}\|Bv  \|_{Z}, 
    \quad \forall v \in \set{v\in \mathfrak{N}_{\theta,n} \cap X  |  Bv  \in Z_{M_b}}.
    \end{split}  
    \end{equation}
    \iffalse 
        \begin{equation}
        \begin{split}
            \|A v -f\|_{Y_{M_r}} &\ge \frac{1}{2}\|Av-f \|_{Y}, 
        \quad 
        \forall v \in  Y_A,
        \\
        \|Bv-g\|_{Z_{M_b}} &\ge \frac{1}{2}\|Bv-g \|_{Z}, 
        \quad \forall v \in Z_B
        \end{split}  
    \end{equation}
    \fi 
   
     }
\end{assu}

{Here $M_r$ and $M_b$ depend on $n$ and may increase with $n$.
}
%\begin{rem}\vskip -15pt
The constant 
$1/2$ can be replaced with any constant larger than $0$, while the constant should not depend on $n$, $M_r$ and $M_b$. %
For example, one can replace $1/2$ with $1-\epsilon$, where  $\epsilon\in (0,1)$ is independent of $n$, $M_r$ and $M_b$. 
%\end{rem}     

 %%%%%%%%%%%%%%%%%%%%%%%%%%%%%%%%%%%%%%%%%%
 
% \newpage 
%{\color{blue}
 \begin{thm}[\textbf{Error estimates of discrete RM} I]\label{thm:convergence-nn-linear-prior-posteriori-generic-discrete}
   Let $Y=L^2(\Omega)$, $Z= L^2(\Gamma )$ and $V = X$.
   Let 
Assumptions~\ref{assum:norm-equivalence-generic} and ~\ref{assum:discrete-norm-lower-bound} be valid. 
{
Let $n, n_1$ be given positive integers.
Let $M_r, M_b$ be chosen according to Assumption~\ref{assum:discrete-norm-lower-bound} with $n_1$.} For $\tau \ge 1$, let $u_{\bN,n}^{\tau,M}$ be a minimizer of $\mathcal{J}_{\tau}^M$ \eqref{def:discrete-pinn-loss}
over the network class $\mathfrak{N}_{\theta,n} \cap X$.
Assume that 
$\norm{v}_{Y_{M_r}}\leq C_3 \norm{v}_Y$ for $v\in Y_{M_r}$ and $\norm{v}_{Z_{M_b}}\leq C_3\norm{v}_Z$ for $v\in Z_{M_b}$, where $C_3>0$ is independent of $M_r$ and $M_b$.
{
Assume further that there exists $\tilde{v} \in X$ such that 
$
\tilde{v} - u_{\mathbb{N},n}^{\tau,M} \in \tilde{V}_{n_1}:=\{v \in \mathfrak{N}_{\theta,n_1} \cap X : Av \in Y_{M_r}, Bv \in Z_{M_b}\}.$} Then, the following error estimates hold
\begin{eqnarray*}
  \norm{u_{\bN,n}^{\tau,M}-u^*}_V  
   &\leq&   2\sqrt{2}C_1^{-1}  \big(\mathcal{J}_\tau^{M}(u_{\bN,n}^{\tau,M}) )^{1/2} + 
  3\sqrt{2}C_1^{-1}  C_3\epsilon_{f,g,u_{\bN,n}^{\tau,M}}^{\frac{1}{2}}.
  %, \text{ and }   \\
  %  \norm{u_{\bN,n}^{\tau,M}-u^*}_V  
  % &\leq&   2\sqrt{2}C_1^{-1} C_3 (1+\tau)^{1/2}  \inf_{v\in \mathfrak{N}_{\theta,n}\cap X}   \big(\mathcal{J}_1 ({v}) )^{1/2} +  \color{blue} 3\sqrt{2}C_1^{-1} C_3\epsilon_{f,g,u_{\bN,n}^{\tau,M}}^{\frac{1}{2}}.
\end{eqnarray*} 
{
 Here $\epsilon_{f,g, u_{\bN,n}^{\tau,M}} = \inf (\norm{A\tilde{v} -f}_Y^2+  \norm{B\tilde{v}-g}_Z^2)$
 where the infimum is taken over
 all $\tilde{v} \in X$
 satisfying $\tilde{v}-u_{\mathbb{N},n}^{\tau,M} \in \tilde{V}_{n_1}$.
}
%{\color{red}
%Here we  denote $\epsilon_{f,g, u_{\bN,n}^{\tau,M}} = \inf_{{v}_{f,g}\in \mathfrak{N}_{\theta,n_1}
%,\,  v_{f,g } -u_{\bN,n}^{\tau,M} \in \mathfrak{N}_{\theta,n_2}}(\norm{Av_{f,g } -f}_Y^2+  \norm{Bv_{f,g}-g}^2)$.
%}
\end{thm}
%}
\begin{proof}  
{ 
Since $\tilde{v} - u_{\bN,n}^{\tau,M} \in \tilde{V}_{n_1}$,
by Assumptions~\ref{assum:norm-equivalence-generic} and~\ref{assum:discrete-norm-lower-bound},   we have 
\begin{eqnarray*}
  C_1 \norm{\tilde{v} - u_{\bN,n}^{\tau,M}}_V  &\leq &  \norm{A(\tilde{v}-u_{\bN,n}^{\tau,M} )}_{Y} + 
  \norm{B(\tilde{v}- u_{\bN,n}^{\tau,M})}_Z\\
   &\leq & 2\norm{A(\tilde{v}-u_{\bN,n}^{\tau,M} )}_{Y_{M_r}} + 
   2\norm{B(\tilde{v}- u_{\bN,n}^{\tau,M})}_{Z_{M_b}}\\
&\leq &2\norm{f - Au_{\bN,n}^{\tau,M}}_{Y_{M_r}} + 2\norm{g- B u_{\bN,n}^{\tau,M}}_{Z_{M_b}}  + 2\sqrt{2}\big(\mathcal{J}_1^{M}(\tilde{v}) )^{1/2} \\  
&\leq &
 2 \sqrt{2}\big(\mathcal{J}_1^{M}(u_{\bN,n}^{\tau,M}) )^{1/2}+ 2\sqrt{2}\big(\mathcal{J}_1^{M}(\tilde{v}) )^{1/2}.
\end{eqnarray*} 
}
Then by the triangle inequality and Assumption \ref{assum:norm-equivalence-generic},  
\begin{eqnarray*}
   \norm{u^*- u_{\bN,n}^{\tau,M}}_V &  \leq &     \norm{u^* -\tilde{v} }_V   +  \norm{\tilde{v} - u_{\bN,n}^{\tau,M}}_V \\
   &\leq&  2 \sqrt{2}C_1^{-1}\big(\mathcal{J}_1^{M}(u_{\bN,n}^{\tau,M}) )^{1/2}+ 3\sqrt{2}C_1^{-1} \big(\mathcal{J}_1^M(\tilde{v}) )^{1/2}
\end{eqnarray*}
By the assumption, we have  $
\mathcal{J}_1^{M}(v) \leq C_3 \mathcal{J}_1(v) $ and  we then obtain the desired conclusion.
\end{proof}

{
\begin{rem}[Existence of $\tilde{v}$ in  Theorem  \ref{thm:convergence-nn-linear-prior-posteriori-generic-discrete}] For Gaussian radial basis networks in Example \ref{exm:bernstein-ineq-guassian-discrete-norm-relation}, 
%	$v_{f,g}$ 
$\tilde{v}$
can be found from the same set  of the  Gaussian radial neural networks for the approximation of the solution $u$ ($n_1=n$). 
For ReLU networks, the existence of $\tilde{v}$
%$v_{f,g}$ 
results from the fact that the summation of ReLU networks is still a ReLU network.
\end{rem}	
}

\begin{exm}[{Gaussian radial neural networks}]\label{exm:bernstein-ineq-guassian-discrete-norm-relation}
Consider the problem of $Au  =f$ on $\Real^d $ with vanishing  $u$ when $\abs{x}\to \infty$.   
Here the operator $A = -\Delta + {\rm Id}$, (Id is the identity operator) and $Bu=0$ when $\abs{x}\to\infty$. 
Then Assumption \ref{assum:discrete-norm-lower-bound} is satisfied for  the following Gaussian radial basis networks
$G_{n,m}(x)$:
\begin{equation*}
    \sum_{k=1}^n a_k \exp(-\abs{x-x_k}^2) :
      a_k \in \mathbb{R},\,
    x_k \in \mathbb{R}^d,\, \inf_{i\neq j}\abs{x_i-x_j}>\frac{1}{m},\, \max_{1\leq k\leq n}\abs{x_k}\leq c m,
\end{equation*} 
 where 
$n,m$ are user-defined  integers, $c>1$ is a constant, and $a_k$ and $x_k$ are unknown.
\end{exm}

	{
To establish  the convergence using the discrete norm relation, we will need  to use a proper $Y_{M_r}$ and the compact set in the relation, which we find below.  

Let $I_N^c$ be the Hermite-Gaussian interpolation operator in 1d ($x\in\Real$) such that 
$I_N^c v  (x ^{(j)}) =  v(x^{(j)})$ for all $v= Q_N(x)\exp(-\frac{\abs{x}^2}{2})$, where $Q_N$ is a polynomial of order no larger than $N$ and $x^{(j)}$'s are the zeros of  the normalized Hermite polynomial $h_{N+1}$.  When $x\in \Real^d$, we still  use the notation $I_N^c$ to represent the  $d$-dimension interpolation operator using the tensor product of one-dimensional interpolation. % where $N$ is the  order of $d$-dimensional polynomials. 
For continuous  $v$, define 
\begin{equation}\label{eq:weighted-interpolation}
\mathcal{I}_N v= \exp (-\frac{\abs{x}^2}{2})I_N^c (\exp (\frac{\abs{x}^2}{2})v).
\end{equation}

\begin{thm} \label{THM:INTP-DERIV-NNs}
	Let
	$\mathcal{I}_N$ be the interpolation operator defined in \eqref{eq:weighted-interpolation}. 
	For any $v$ in the form of $ G_{n,m}$, there exist constants $c,c_1,c_2>0$ independent of $n$ and $m$ and $c_2 \gg 1$ such that 
	\[ \norm{A v -  \mathcal{I}_N Av }_{L^2}\leq c\exp (-c_1m^2)\norm{v}_{L^2}, \qquad N = c_2m^2.\]
\end{thm} 
\begin{proof}
	The proof can be found in \ref{APP:THM:INTP-DERIV-NNs}.
\end{proof}

Let $v$ be in the form of  $G_{n,m}$. By  Theorem \ref{THM:INTP-DERIV-NNs}, we can find $Y_{M_r}$ by 
observing that  
\[\int_{\Real^d}(\mathcal{I}_N Av )^2\,dx =\int_{\Real^d}({I}_N^c \exp(\frac{\abs{x}^2}{2}) Av ) )^2\exp(-\abs{x})^2\,dx.\]
By the definition of $I_N^c$, we can find a quadrature rule 
$\{(x^{(i)},w_i)\}_{i=1}^{M_r}$ ($M_r=O(N^d)$, for example) such that 
\begin{equation}
\int_{\Real^d}(\mathcal{I}_N Av )^2\,dx = 
\sum_{i=1}^{M_r} \left(Av(x^{(i)}) \right)^2w_i.
\end{equation}  
Let $Y=L^2(\mathbb{R}^d)$ and $Y_{M_r}\subset Y$ be  equipped with the  discrete norm $\norm{v}_{Y_{M_r}}= \left(\sum_{i=1}^{M_r} \big(v(x^{(i)})\big)^2w_i \right)^{1/2}$.  
Then by the Cauchy-Schwarz inequality and Theorem \ref{THM:INTP-DERIV-NNs}, we have, for $v$ in the form of $G_{n,m}$,
\begin{eqnarray*}
\abs{ \norm{Av}_Y^2-  \norm{Av}_{Y_{M_r}} ^2} &=&
\abs{\norm{Av}_Y^2-\norm{\mathcal{I}_N Av}_{Y}^2}\\
&\leq&  
   \norm{Av -  \mathcal{I}_N Av}_{L^2}\norm{Av + \mathcal{I}_N Av}_{L^2}    \leq C \exp (-c_1m^2)\norm{v}_{L^2}\norm{Av}_{L^2}, 
\end{eqnarray*}
where $C$ is a constant depending on $c$ from Theorem \ref{THM:INTP-DERIV-NNs}.
It can be checked by energy estimates that $\norm{v}_{L^2}\leq \norm{Av}_{L^2}$. 
Then we have 
$\abs{ \norm{Av}_Y^2-  \norm{Av}_{Y_{M_r}} ^2} \leq C \exp(-c_1 m^2)\norm{Av}_{L^2}^2$. 
Picking $m$ such that  $C \exp(-c_1 m^2)\leq \frac{3}{4} $, we obtain  that
$ 4\norm{Av}_{Y_{M_r}} ^2\geq \norm{Av}_{Y=L^2}^2\geq 4/7\norm{Av}_{Y_{M_r}} ^2$.
 
Then for any Gaussian neural networks $v\in W^{2,2}$ such that $Av \in Y_{M_r}$, the first inequality with $f=0$ in  Assumption \ref{assum:discrete-norm-lower-bound} is satisfied. 
Let $v=\sum_{k=1}^n a_k \exp(-\abs{x-x_k}^2) $ and
let $v_f=\sum_{k=1}^n f_k \exp(-\abs{x-x_k}^2)$ be an approximation of $f$. 
%and let $v=\sum_{k=1}^n a_k \exp(-\abs{x-x_k}^2)\in G_{n,m}$. 
By the above discussion, we have the following discrete norm relation: 
\begin{equation*}
 \frac{2}{\sqrt{7}} \norm{A(v-v_f)}_{Y_{M_r}}\leq  \norm{A(v-v_f)}_{L^2} \leq 2\norm{A(v-v_f)}_{Y_{M_r}}.
\end{equation*}
Then the error estimates of Gaussian radial  basis networks for the problem $Av=f$ and $Bu=0$ when $\abs{x}\to\infty$ can derived by Theorem \ref{thm:convergence-nn-linear-prior-posteriori-generic-discrete}.
} 

%\begin{rem}
In general,  Assumption \ref{assum:discrete-norm-lower-bound} is not readily verified using the above approach via the inverse estimate (Bernstein-type inequality). In fact, the Bernstein-type inequality for deep feed-forward neural networks is unavailable for even simple ReLU networks (see a counterexample in \cite{Hong2021priori}) when there are no constraints on the weights and biases.
In the next subsection, we will use the Rademacher complexity 
to capture the effect of the discretization.
%of a set of networks for quasi-minimizers. 
The Rademacher complexity for two-layer networks may suggest a lift of the curse of dimensionality as it depends logarithmically on the dimension $d$. 
%\end{rem}
%%-----------------------------------------

%%%%%%%%%%%%%%%%%%%%%%%%%%%%%%%%%%%%%%%%%%%%%%%%%%%%%%%%%%%

%% file: drm.tex
%%%%%%%%%%%%%%%%%%%%%%%%%%%%%%%%%%%%%%%%%%%%%%%%%%%%%%%%%%%
 
\subsection{Using Rademacher complexity} \label{sec:prototype-PINNs}

The analysis presented in this section will be based on
the uniform law of large numbers
for well-designed classes of neural networks.

\begin{defn}[Rademacher complexity, \cite{gnecco2008approximation,Wainwright2019high}]
    Given a collection $\{X_i\}_{i=1}^M$ of i.i.d. random samples, 
    the Rademacher complexity of the function class
    $\mathfrak{F}$ is defined by
   $ % \begin{equation}
        R_M(\mathfrak{F}) = \mathbb{E}_{\{X_i,\epsilon_i\}_{i=1}^M}\left[ \sup_{f \in \mathfrak{F}} |\frac{1}{M}\sum_{i=1}^M \epsilon_i f(X_i)| \right],
  $ %   \end{equation}
    where $\epsilon_i$'s are i.i.d. Rademacher random variables 
    i.e., $\bP(\epsilon_i = 1) = \bP(\epsilon_i = -1) = 0.5$.
\end{defn}

For many function classes,  the upper bounds of the Rademacher complexity  are known.
For example, 
see \cite{Neyshabur2018towards} for
a class of two-layer neural networks
and 
see \cite{ma2019barron}
for a class of Barron functions.
In many cases, it can be shown that 
the Rademacher complexity 
converges to zero as the number of samples $M$
grows to infinite, such as in \cite{ma2019barron} where the convergence rate is half.

Let $\{u_k^*\}$ be a sequence approximating the solution, as  defined in 
Definition~\ref{def:solution}.
By definition, the sequences $\{Au_k^* - f\}$ and $\{Bu_k^* - g\}$
are uniformly bounded.
In what follows, we make an additional assumption on these sequences required for the uniform law of large numbers.

\begin{assu} \label{assum:Rademacher}
    Let $u^*$ be a solution to \eqref{eq:linear-generic-bvp}
    in the sense of Definition~\ref{def:solution}
    and $\{u^*_k\}$ be its corresponding sequence in $X$.
    Let 
    {$G_r \ge \max\{ \|f\|_{L^\infty(\Omega)}, \sup_{k} \|Au_{k}^* - f\|_{L^\infty(\Omega)}\},
    $ and $
    G_b \ge \max\{\|g\|_{L^\infty(\Gamma)}, \sup_{k}
    \|Bu_k^* - g\|_{L^\infty(\Gamma)}\}$.
	Assume $G_r, G_b < \infty$.}
\end{assu}

Under the Assumption~\ref{assum:Rademacher},
we introduce the following function classes.
\begin{defn} \label{def:ULLN-function-class}
    Let $G_r, G_b$ be positive numbers defined in 
    Assumption~\ref{assum:Rademacher}.
    Define a subclass of $\mathfrak{N}_{\theta,n}$ by 
  	{$\mathfrak{N}_{\theta,n}^{\mathbb{Q}} = \{\mathcal{R}[\theta] \in \mathfrak{N}_{\theta,n} : \theta \in \Theta_{\mathbb{Q}}(\vec{\bm{n}}_n) \}$},
    where $\Theta_{\mathbb{Q}}(\vec{\bm{n}}_n)$
    is similarly defined 
    as in \eqref{def:NN-parameter-space}
    but all weights and biases are rational numbers. 
    Introduce two function classes: \vskip -15pt 
    \begin{eqnarray*} 
            \mathfrak{F}_{r,n} &:=& \{ Av - f \in Y : v \in \mathfrak{N}_{\theta,n}^{\mathbb{Q}}\cap X \text{ and } 
            \|Av - f\|_{L^\infty(\Omega)} \le G_r
            \}, 
            \\
            \mathfrak{F}_{b,n} &:=& \{ Bv - g \in Z : v \in \mathfrak{N}_{\theta,n}^{\mathbb{Q}} \cap X
            \text{ and } 
            \|Bv - g\|_{L^\infty(\Gamma)} \le G_b
            \}. 
    \end{eqnarray*} 
% \vskip -12pt  \noindent
    Finally, define a subclass of $\mathfrak{N}_{\theta,n}^\mathbb{Q}$
    \begin{equation}
        \widetilde{\mathfrak{N}}_{\theta,n}^{\mathbb{Q}}
        = \{v \in \mathfrak{N}_{\theta,n}^{\mathbb{Q}} \cap X :
            Av - f \in \mathfrak{F}_{r,n}, 
            Bv- g \in \mathfrak{F}_{b,n} \}.
    \end{equation}
\end{defn}

From Assumption~\ref{assum:Rademacher}, the countable set
$\widetilde{\mathfrak{N}}_{\theta,n}^{\mathbb{Q}}$ is not empty
{ as it contains the zero function.} 
One can also consider an uncountable function class
by strengthening Assumption~\ref{assum:Rademacher}
with the uniform topology $C^0$.

Next, we  bound the discrepancy between the continuous norm and its discretization
using the Rademacher complexity.
The following lemma is obtained by applying the uniform law of large numbers 
on the function class defined in Definition~\ref{def:ULLN-function-class}.

\begin{lem} \label{lemma:Rademarcher-Loss}
	 
    Suppose Assumption~\ref{assum:Rademacher}
    holds.
    Let $G_r, G_b$ be positive numbers defined in 
    Assumption~\ref{assum:Rademacher}.
    Let $\{x_i^r\}_{i=1}^{M_r}$
    and $\{x_i^b\}_{i=1}^{M_b}$
    be i.i.d. samples 
    following probability densities $\rho$ over $\Omega$ and $\rho_b$ over $\Gamma$,
    respectively in the discrete RM loss \eqref{def:discrete-pinn-loss}.
    Let $Y = L^p_{\rho}(\Omega)$ and $Z = L^p_{\rho_b}(\Gamma)$ for $p \ge 1$.
    For any small $\delta_r, \delta_b > 0$, 
    \begin{align}
    	\sup_{Av-f \in \mathfrak{F}_{r,n}} \left|\|Av-f\|_{Y_{M_r}}^p - \|Av-f\|_{L_{\rho}^p(\Omega)}^p \right|
    	&\le 2R_{M_r}(\mathfrak{F}_{r,n}^p) + \frac{\delta_r}{2}, 
    	\text{ with prob. at least } Q_r, \label{eq:rad-complex-r}\\
    	\sup_{Bv-g \in \mathfrak{F}_{b,n}} \left|\|Bv-g\|_{Z_{M_b}}^p - \|Bv-g\|_{L_{\rho_b}^p(\Gamma)}^p \right|
    	&\le 2R_{M_b}(\mathfrak{F}_{b,n}^p) + \frac{\delta_b}{2},   
    	\text{ with prob. at least } Q_b, \label{eq:rad-complex-b}
    \end{align}
	where $Q_r=(1-2\exp(-\frac{M_r\delta_r^2}{32G_r^{2p}}))$,  $Q_b=(1-2\exp(-\frac{M_b\delta_b^2}{32G_b^{2p}}))$,
	$\mathfrak{F}_{r,n}^p = \{v^p : v \in \mathfrak{F}_{r,n}\}$,
	and $\mathfrak{F}_{b,n}^p = \{v^p : v \in \mathfrak{F}_{b,n}\}$.
	Furthermore, with probability at least $Q_{r,b}=Q_r Q_b$, 
    \begin{align*}
        \sup_{v \in \widetilde{\mathfrak{N}}_{\theta,n}^{\mathbb{Q}}}  \left|\mathcal{J}_\tau^M(v) - \mathcal{J}_\tau(v) \right|
        \le
        2R_{M_r}(\mathfrak{F}_{r,n}^p) 
        +
        2\tau R_{M_b}(\mathfrak{F}_{b,n}^p) + \frac{\delta_r}{2} + \tau \frac{\delta_b}{2}.
    \end{align*}
\end{lem}
\begin{proof}
    Recall from Section~\ref{sec:ssec-effect-discretizing-integrals} that 
    \begin{align*}
        \|Av-f\|_{Y_{M_r}}^p = \frac{1}{M_r}\sum_{i=1}^{M_r} (f(x_i^r) - Av(x_i^r))^p,  \quad 
        \|Bv-g\|_{Z_{M_b}}^p = \frac{1}{M_b}\sum_{i=1}^{M_b} (g(x_i^b) - Bv(x_i^b))^p.
    \end{align*}
    We observe that 
    \begin{eqnarray*}
         \left|\mathcal{J}_\tau^M(v) - \mathcal{J}_\tau(v) \right|
         \le
        \left|\|Av-f\|_{Y_{M_r}}^p - \|Av-f\|_{L_{\rho}^p(\Omega)}^p \right|
        +
        \tau \left|\|Bv-g\|_{Z_{M_b}}^p - \|Bv-g\|_{L_{\rho_b}^p(\Gamma)}^p \right|.
    \end{eqnarray*}
    For each $F \in \mathfrak{F}_{r,n}$, 
    since $\|F\|_{L^\infty(\Omega)} \le G_r$,
    $\bP(|F(x)|\le G_r) = 1$.
    Since $\mathfrak{F}_{r,n}$ is countable,   
    we have 
    $\bP(\sup_{F \in \mathfrak{F}_{r,n}} |F(x)|\le G_r) = 1$ by the continuity of the probability measure.
    A similar argument leads to the conclusion for $\mathfrak{F}_{b,n}$.
    Then \eqref{eq:rad-complex-r} and \eqref{eq:rad-complex-b} can be obtained by invoking a uniform law via the Rademacher complexity (e.g. Theorem 4.2 of \cite{Wainwright2019high}).  
    %with probability at least $1-2\exp(-\frac{M_b\delta_b^2}{32G_b^{2p}})$.
    Combining the above two estimates leads to the last desired conclusion.
\end{proof}

{
\begin{exm}[Rademacher complexity for two-layer networks]\label{exm:Rademacher}
Let $\mathfrak{N}_{\mathbb{N},n}$ be the class of two-layer neural networks defined by
	$$
	\mathfrak{N}_{\mathbb{N},n} := \left\{ [0,1]^d \ni x \mapsto c_0 + \sum_{i=1}^n c_i\phi(w_i^\top x + b_i) : \max |b_i|, \|w_i\|_{\ell_1}, |c_i| \le \omega_{\max} \right\},
	$$
	where $\phi$ is $\gamma$-Lipschitz and anti-symmetric, i.e., $\phi(-x)=-\phi(x)$ (e.g. $\phi(x) = \tanh(x)$).
	For the sake of the length of the paper, we briefly give a sketch of how one can estimate $R_{M_r}(\mathfrak{F}_{r,n}^p)$ 
	and
	$R_{M_b}(\mathfrak{F}_{b,n}^p)$ 
	when
	$A$ is the Laplacian $\Delta$ on $\Omega=[0,1]^d$ and $B$ is the identity operator on $\partial\Omega$.
	It then can be checked that 
        (e.g. \cite{neyshabur2018the})
	%\cite{rebeschini20}
%	https://www.stats.ox.ac.uk/~rebeschi/teaching/AFoL/20/material/lecture03.pdf
	\begin{align*}
		R_{M}(\mathfrak{N}_{\mathbb{N},n})
		\le \frac{\omega_{\max}}{\sqrt{M}}\left(1 + 2\gamma \omega_{\max} (1 + \sqrt{2\log (2d)}) \right). 
	\end{align*}
	Observe that 
	$\mathfrak{F}_{r,n}$  and $\mathfrak{F}_{b,n}$ are given by
	\begin{align*}
		\mathfrak{F}_{r,n} &=\left\{ [0,1]^d \ni x \mapsto 
		\sum_{i=1}^n c_i\|w_i\|_{\ell_2}^2\phi''(w_i^\top x + b_i) - f(x) : \max |b_i|, \|w_i\|_{\ell_2}, |c_i| \le \omega_{\max} \right\}, \\
		\mathfrak{F}_{b,n} &= \left\{ [0,1]^d \ni x \mapsto c_0 + \sum_{i=1}^n c_i\phi(w_i^\top x + b_i) - g(x) : \max |b_i|, \|w_i\|_{\ell_1}, |c_i| \le \omega_{\max} \right\}.
	\end{align*}
	Since $\|g\|_{L^\infty(\Gamma)} < \infty$, 
	$R_{M_b}(\mathfrak{F}_{b,n})$ is bounded by
	$R_{M_b}(\mathfrak{N}_{\mathbb{N},n}) + \frac{\|g\|_{L^\infty(\Gamma)}}{\sqrt{M_b}}$
	(Chapter 4 of \cite{Wainwright2019high}).
	Also, since $\|f\|_{L^\infty(\Omega)}  < \infty$,
	$R_{M_r}(\mathfrak{F}_{r,n})$ is bounded by
	$R_{M_r}(\Delta \mathfrak{N}_{\mathbb{N},n}) + \frac{\|f\|_{L^\infty(\Omega)}}{\sqrt{M_r}}$
	where
	$\Delta \mathfrak{N}_{\mathbb{N},n} = \{ \Delta u : u \in \mathfrak{N}_{\mathbb{N},n}\}$.
	Assuming $\phi''$ is $\gamma'$-Lipschitz and anti-symmetric,
	$R_{M_r}(\Delta \mathfrak{N}_{\mathbb{N},n})$ can also be estimated similarly.
	Lastly, since every function in $\mathfrak{F}_{r,n}$ (or  $\mathfrak{F}_{b,n}$)
	is bounded and $x\mapsto x^p$ is Lipschitz on a bounded interval,
	an upper bound of the Rademacher complexity of $\mathfrak{F}_{r,n}^p$ (or $\mathfrak{F}_{b,n}^p$) can be found.
\end{exm}
}

With Lemma~\ref{lemma:Rademarcher-Loss} 
and Theorem~\ref{thm:convergence-nn-linear-prior-posteriori-generic},
we establish error estimates for the discrete RM. %
\begin{thm}[\textbf{Error estimates for discrete RM} II] \label{thm:err-discrete-pinns}
Suppose Assumptions
\ref{assum:solution-existence},
\ref{assum:norm-equivalence-generic}
and
\ref{assum:Rademacher}
hold. 
Suppose $\{x_i^r\}_{i=1}^{M_r}$ and $\{x_i^b\}_{i=1}^{M_b}$
are be i.i.d. samples 
following probability densities $\rho$ over $\Omega$ and $\rho_b$ over $\Gamma$,
respectively in the discrete RM loss \eqref{def:discrete-pinn-loss} where $\tau \ge 1$.
Let $Y = L^p_{\rho}(\Omega)$ and $Z = L^p_{\rho_b}(\Gamma)$ for $p \ge 1$.
Let $u_{\mathbb{N},n}^{\tau,M}$ be 
a solution to 
$ %\begin{equation*}
    \min_{u \in \widetilde{\mathfrak{N}}_{\theta,n}^{\mathbb{Q}}} \mathcal{J}_{\tau}^{M}(u),
$ %\end{equation*}
where $M=(M_r,M_b)$.
Let $u^*$ be the solution to \eqref{eq:linear-generic-bvp}
in the sense of Definition~\ref{def:solution}.
%and $\{u^*_k\}$ be its corresponding sequence in $X$.
Then, for $\delta > 0$,
with probability at least $Q_{r,b}=(1-2\exp(-\frac{M_r\delta^2}{32G_r^{2p}}))(1-2\exp(-\frac{M_b\delta^2}{32G_b^{2p}}))$, 
we have
\begin{align*}
        \norm{u_{\mathbb{N},n}^{\tau,M} - u^*}_{V} 
        \leq  C_1^{-1}2^{\frac{p-1}{p}} \left[\mathcal{J}_{\tau}^M(u_{\mathbb{N},n}^{\tau,M}) + 2\tilde{R}_M(\mathfrak{F}_n^p) + (1+\tau)\delta/2 \right]^{\frac{1}{p}},
\end{align*}
where $\tilde{R}_M(\mathfrak{F}_n^p) := {R}_{M_r}(\mathfrak{F}_{r,n}^p) + \tau  {R}_{M_b}(\mathfrak{F}_{b,n}^p)$.
Also, for any $\epsilon > 0$, there exists $u^*_{\epsilon} \in X$
such that with probability  $Q_{r,b}$ (at least), 
we have
\begin{eqnarray*}
            \norm{u_{\mathbb{N},n}^{\tau,M} - u^*}_{V} 
         \le C_1^{-1}2^{\frac{p-1}{p}} \left[(1+\tau)\inf_{w \in \tilde{\mathfrak{N}}_{\theta,n}^{\mathbb{Q}}}
        (C_2\norm{w - u_\varepsilon^*}_X
        + \epsilon)^p + 4\tilde{R}_M(\mathfrak{F}_n^p) + (1+\tau)\delta \right]^{\frac{1}{p}}. 
\end{eqnarray*}
Here $C_1$ and $C_2$ are the constants defined in 
\eqref{eq:norm-equivalence-generic}.
\end{thm}
\begin{proof}
    From Theorem~\ref{thm:convergence-nn-linear-prior-posteriori-generic},
    we have 
   $ % \begin{align*}
        \norm{u_{\mathbb{N},n}^{\tau,M} - u^*}_{V}\leq  C_1^{-1}2^{\frac{p-1}{p}} \left[\mathcal{J}_{\tau}(u_{\mathbb{N},n}^{\tau,M})\right]^{\frac{1}{p}}.
    $ % \end{align*}
    Since $u_{\mathbb{N},n}^{\tau,M} \in \widetilde{\mathfrak{N}}_{\theta,n}^{\mathbb{Q}}$,
    it follows from Lemma~\ref{lemma:Rademarcher-Loss} that 
    with probability  $Q_{r,b}$ (at least), 
    we have
    \begin{align*}
        \norm{u_{\mathbb{N},n}^{\tau,M} - u^*}_{V} 
        \leq  C_1^{-1}2^{\frac{p-1}{p}} \left[\mathcal{J}_{\tau}^M(u_{\mathbb{N},n}^{\tau,M}) + 2\tilde{R}_M(\mathfrak{F}_n^p) + (1+\tau)\delta/2 \right]^{\frac{1}{p}}.
    \end{align*}
%    where $\tilde{R}_M(\mathfrak{F}_n) = {R}_{M_r}(\mathfrak{F}_{r,n}) + {R}_{M_b}(\mathfrak{F}_{b,n})$.
%    Since %$u_{\bN,n} ^{\tau,M}$ is a minimizer of the loss and 
%    $\mathcal{J}_{\tau}^M(w)\leq \mathcal{J}_{\tau}(w)$, we deduce that 
   %
   Similarly, with probability  $Q_{r,b}$ (at least), 
   we have
    \begin{equation} \label{eq:bound-opt-discretePINNs}
        \begin{split}
            \norm{u_{\mathbb{N},n}^{\tau,M} - u^*}_{V} 
        %&\le  C_1^{-1}2^{\frac{p-1}{p}} \left[\inf_{w \in \tilde{\mathfrak{N}}_{\theta,n}^{\mathbb{Q}}}\mathcal{J}_{\tau}^M(w) + 2\tilde{R}_M(\mathfrak{F}_n) + (1+\tau)\delta/2 \right]^{\frac{1}{p}} \\
        &\le C_1^{-1}2^{\frac{p-1}{p}} \left[\inf_{w \in \tilde{\mathfrak{N}}_{\theta,n}^{\mathbb{Q}}}\mathcal{J}_{\tau}(w) + 4\tilde{R}_M(\mathfrak{F}_n^p) + (1+\tau)\delta \right]^{\frac{1}{p}}.
        \end{split}
    \end{equation}
    For $\epsilon > 0$,
    let $u_{\epsilon}^* \in X$ be an approximation of $u^*$ such that
    $\|Au_{\epsilon}^* - f\|_Y + \|Bu_{\epsilon}^*-g\|_Z < \epsilon$. 
    % Let $ u_\varepsilon^*$ be a mollified $u^* \in X$
    % satisfying 
    % $\norm{ Au_\varepsilon^* -A u^*}_Y +  \norm{B u_\varepsilon -u^*}_Z < \epsilon$
   % (Assumption~\ref{assum:solution-existence}).
%
    From \eqref{eq:bound-optimization} and \eqref{eq:tri-error-mollif-appr}, we have 
    \begin{align*}
	\mathcal{J}_{\tau}(w\big)
	\le  (1+\tau)(\norm{Aw - f}_Y 
	+\norm{Bw -g}_Z)^p
	\le (1+\tau)(C_2\norm{w - u_\varepsilon^*}_X
 + \epsilon)^p.
    \end{align*}

Combining the above inequality with \eqref{eq:bound-opt-discretePINNs}, we obtain
\begin{eqnarray*}
            \norm{u_{\mathbb{N},n}^{\tau,M} - u^*}_{V} 
        &\le C_1^{-1}2^{\frac{p-1}{p}} \left[(1+\tau)\inf_{w \in \tilde{\mathfrak{N}}_{\theta,n}^{\mathbb{Q}}}
        (C_2\norm{w - u_\varepsilon^*}_X
        + \epsilon)^p + 4\tilde{R}_M(\mathfrak{F}^p_n) + (1+\tau)\delta \right]^{\frac{1}{p}}.
 \end{eqnarray*}
The proof is then complete.
\end{proof}

Next, we characterize
the conditions under which
a sequence of minimizers of the discrete RM loss functionals  
converges to the solution 
strongly in $V$.
\begin{thm}[\textbf{Convergence of discrete RM}] \label{thm:convg-discretePINNs}
    Under the same conditions of Theorem~\ref{thm:err-discrete-pinns},
    suppose 
    $\lim_{M_r \to \infty} {R}_{M_r}(\mathfrak{F}_{r,n}^p) = 0$,
    and $\lim_{M_b \to \infty} {R}_{M_b}(\mathfrak{F}_{b,n}^p) = 0$
    for all $n$.
    {
    Suppose further that 
    there exists a solution sequence $\{v_k^*\}$
    (Definition~\ref{def:solution})
    that belongs to 
    $\overline{\bigcup_{n=1}^\infty \widetilde{\mathfrak{N}}^{\mathbb{Q}}_{\theta,n}}$ in the topology of $(X,\|\cdot\|_X)$.}
    Then,
    we have
    \begin{align*}
        \lim_{n \to \infty } \lim_{M \to \infty} 
        \norm{u_{\mathbb{N},n}^{\tau, M} - u^*}_{V} 
        = 0, \quad M=(M_r, M_b),
    \end{align*} 
	in probability.
\end{thm}
\begin{proof}
	{
	Let $\{v_k^*\}$ be the solution sequence
	and let $\{\epsilon_k\}$ be a positive decreasing sequence converging to 0.
    For $k \gg 1$, it follows from
    the assumption that
    there exists an integer $n_k$ and a neural network $u^{\mathbb{Q}}_{n_k} \in \widetilde{\mathfrak{N}}^{\mathbb{Q}}_{\theta,n_k} \cap X$ such that 
    $\|u^{\mathbb{Q}}_{n_k}-  v_k^*\|_X \le \epsilon_k$.
    From Theorem~\ref{thm:err-discrete-pinns},
    by choosing $\delta_r = 2M_r^{-\frac{1}{2} + \epsilon}$ and $\delta_b = 2M_b^{-\frac{1}{2} + \epsilon}$
    for $0 < \epsilon < \frac{1}{2}$,
    with probability at least $(1-2\exp(-\frac{M_r^{\epsilon}}{8G_r^{2p}}))(1-2\exp(-\frac{M_b^{\epsilon}}{8G_b^{2p}}))$,
    we have
    \begin{align*}
        \norm{u_{\mathbb{N},n_k}^{\tau,M} - u^*}_{V} 
        \leq  C_1^{-1}2^{\frac{p-1}{p}} \left[\mathcal{J}_{\tau}^M(u^{\mathbb{Q}}_{n_k}) + 2\tilde{R}_M(\mathfrak{F}_n^p) + M_r^{-\frac{1}{2} + \epsilon} + \tau M_b^{-\frac{1}{2} + \epsilon} \right]^{\frac{1}{p}}.
    \end{align*}
    By letting $M_r, M_b \to \infty$, 
    we have
    $\lim_{M\to \infty} \norm{u_{\mathbb{N},n_k}^{\tau,M} - u^*}_{V} 
    \leq  C_1^{-1}2^{\frac{p-1}{p}} \left[\mathcal{J}_{\tau}(u^{\mathbb{Q}}_{n_k}) \right]^{\frac{1}{p}}$
    with probability 1 over i.i.d. samples.
%    \begin{align*}
%        \lim_{M\to \infty} \norm{u_{\mathbb{N},n_k}^{\tau,M} - u^*}_{V} 
%        \leq  C_1^{-1}2^{\frac{p-1}{p}} \left[\mathcal{J}_{\tau}(u^{\mathbb{Q}}_{n_k}) \right]^{\frac{1}{p}}.
%    \end{align*}
    Since 
    \begin{align*}
    	\mathcal{J}_{\tau}(u^{\mathbb{Q}}_{n_k}) &\le (\|f - Av_k^*\|_Y + \|A(v_k^* - u^{\mathbb{Q}}_{n_k})\|_Y )^p + \tau (\|g - Bv_k^*\|_Z +\|B(v_k^* - u^{\mathbb{Q}}_{n_k})\|_Z )^p \\
    	&\le (\|f-Av_k^*\|_Y + C_2\epsilon_k)^p + (\|g-Bv_k^*\|_Z+C_2\epsilon_k)^p,
    \end{align*}
    we have  $\lim_{k\to \infty} \mathcal{J}_{\tau}(u^{\mathbb{Q}}_{n_k}) = 0$.
    Since $\widetilde{\mathfrak{N}}^{\mathbb{Q}}_{\theta,m_1} \subset \widetilde{\mathfrak{N}}^{\mathbb{Q}}_{\theta,m_2}$
	for any $m_1 < m_2$, the proof is completed.}
\end{proof}

 

%% file: hpvpinn.tex
%%===================================================
\section{hp-Variational Residual Minimization} 
\label{sec:hp-vpinn}
%%%%%%%%%%%%%%%%%%%%%%%%%%%%%%%%%%%%%%%%%%%%%%
%The residual minimization formulation discussed in previous sections
%does not perform integration-by-parts. 
%Thus it requires some regularity of the input data and the neural network (activation function). For example, cannot be used for an elliptic problem as it is not twice differentiable.   
 
 In this section, we consider error estimates of the hp-variational residual minimization (hp-VRM), which can accommodate non-smooth networks   for high-order differential operators  by performing integration-by-parts.  
 %
%%%%%%%%%%%%%%%%%%%%%%%%%%%%%%%%%%%
\subsection{Definition of hp-VRM}
  Let $\set{\Omega_k}_{k=1}^\infty$ be
  a partition of the domain $\Omega$.  
  Suppose $\set{ \Phi_{k,i}}_{k,i=1}^\infty$ forms a complete  orthonormal basis in a Hilbert space $Y$.
  Also, let $\set{\Phi_{k,i}}_{i=1}^\infty$ be a complete  orthonormal basis in $Y|_{\Omega_k}$, which is defined with the same structure as in $Y$ but over the domain $\Omega_k$. 
  Then, the loss functional 
  \eqref{eq:loss-functional}
  can be  written as 
%  \begin{equation*}
  $\mathcal{J}_\tau (v) =  \sum_{k,i=1}^{\infty}  (f- A v, \Phi_{k,i})_Y  ^2+  \tau \norm{Bu-g}^2_{Z}$,
  %\end{equation*} 
  and its corresponding truncation is given by 
  \vskip -10pt 
  \begin{equation} \label{eq:loss-function-hp-vpinn}
      \mathcal{J}_\tau^{h,\bm{N}}(v) = 
      \sum_{k=1}^{K_h} \sum_{i=1}^{N_k}
      (f- A v, \Phi_{k,i})_{Y}^2
      +  \tau \norm{Bu-g}^2_{Z},
  \end{equation} 
  \vskip -8pt  \noindent
  where 
  $\bm{N}=(N_1,\cdots, N_{K_h})$.
  For simplicity, 
  we  write $K_h$ as $K$.
  If $N_k = N$ for all $k$, we write $\mathcal{J}_\tau^{h,\bm{N}}$
as $\mathcal{J}_\tau^{h,N}$.
%   In practice, we need to truncate the infinite series: 
% 	\begin{equation}
% \mathcal{J}_{\tau}^{h, N } (v) = \sum_{k=1}^{K_h}\sum_{i=1}^N  (f-A v, \Phi_{k,i})_Y ^2 +\tau    \norm{Bv-g}^2_{Z}.
% \end{equation} 
%
We can   perform integration-by-parts
on  $(f- A v, \Phi_{k,i})_Y$.
The goal of hp-VRM is to find 
a solution to the   minimization problem  
%\begin{equation*}
$\inf_{v \in S} \mathcal{J}_\tau^{h,\bm{N}}(v)$, 
%\end{equation*}
where the feasible space $S$ is to be determined shortly.

When $f -A v$ is merely in $L^p(\Omega)$ ($p\in [1,\infty))$, we cannot use the RM formulation. In contrast, we can still use the hp-VRM formulation by performing integration-by-parts.
In other words, the hp-VRM formulation with the piecewise constant basis is a \textit{weak} formulation of RM.

 \begin{rem} 
 In practice, it is convenient to consider the following functional  
 	%\begin{equation*}%\label{eq:loss-function-hp-vpinn}
 	$\mathcal{J}_{\tau}^{h,N} (v) = \sum_{k=1}^{K}\sum_{i=1}^N\alpha_{k,i}^2  (f- A v, \Phi_{k,i})_Y ^2 +\tau    \norm{Bv-g}^2_{Z},
 	$ %\end{equation*} 
 	where $0<M_0\leq \alpha_{k,i}^2\leq M$ as   
 	$
 	M_0J_0^{h,N} (v) \leq  \sum_{k=1}^{K}\sum_{i=1}^N\alpha_{k,i}^2\ (f- A v, \Phi_{k,i})_Y ^2 \leq M J_0^{h,N}(v).$
 	One has the freedom to choose $\alpha_{k,i}$ and $\Phi_{k,i}$. %Here $\alpha_{k,i}$
 	To simplify the analysis, we always assume that 
 	for each $k$, 
 	$\set{ \Phi_{k,i}}_{i \ge 1}$ is a complete orthonormal basis
 	of $Y|_{\Omega_k}$ and $\alpha_{k,i}=1$.
 	%\textcolor{blue}{This implies $\alpha_{k,i}^2 = \|\Phi_{k,i}\|_Y^{-2}$.}
 \end{rem} 
 
 \begin{rem}
 We can also treat the boundary residual in  a similar way. In practice, the boundary might be irregular. It is then practical to use piecewise constants as the basis, instead of piecewise polynomials.  
 \end{rem}

To illustrate the hp-VRM formulation,
we consider two special cases when $Y= L^2(\Omega)$
and $Y|_{\Omega_k}= L^2( \Omega_k)$. 
The first one is the basis of orthonormal polynomials over smooth regular domains. 
This formulation with Jacobi-type polynomials as the basis  is used in \cite{kharazmi2019variational} for tensor product domains. 

The second one is the basis of piecewise constants (i.e., $N=1$) over possibly very complicated domains, 
leading to a weak formulation of residual minimization. 
The first term in the functional $J_\tau^{h,1}(v)$ becomes 
%\begin{equation}
%    \sum_{k=1}^K \int_{\Omega_k} \frac{1}{|\Omega_k|} (f- Av)^2 \,dx  = 
 $   \sum_{k=1}^K |\langle f-Av, \Phi_{k,1}\rangle_{L^2(\Omega_k)}|^2
    =
    \sum_{k=1}^K \frac{1}{\abs{\Omega_k}}\big(\int_{\Omega_k}  f- Av  \,dx \big)^2.
$
%\end{equation}
Specifically, we employ $\set{\Phi_{k,1}(x)= \abs{\Omega_k}^{-1/2}\mathbb{I}_{\Omega_k}(x)}_{k=1}^K$,
which is an orthonormal basis in $L^2(\Omega)$.
Here $\abs{\Omega_k}$ represents the Lebesgue measure (volume) of the domain $\Omega_k$.  
Let $f$, $v$, and the coefficients of $A$ be smooth enough so that $f-Av$ is continuous (or in $L^\infty$).
Then, by the mean value theorem of definite integrals, there exists a point $x_k^* \in \Omega_k $ such that the above term becomes 
%\begin{equation} 
% \sum_{k=1}^K \frac{1}{\abs{\Omega_k}}\big(\int_{\Omega_k} (f(x_k^*)- A[v](x_k^*))^2,
$\sum_{k=1}^K \frac{1}{\abs{\Omega_k}}(f(x_k^*)- A[v](x_k^*))^2$,
%\end{equation} 
which is the first term  in the discrete RM formulation \eqref{def:discrete-pinn-loss}. 
If 
$f-Av$ is not in $L^\infty$, we may apply integration-by-parts 
in the following formula:
%\begin{equation}
 $   \sum_{k=1}^K \frac{1}{\abs{\Omega_k}}\big(\int_{\Omega_k} (f- Av) \,dx \big)^2 +   \tau\sum_{l=1}^M \frac{1}{\abs{\Gamma_l}}\big(\int_{\Gamma_l} (Bu- g) \,dx \big)^2$.
%\end{equation} 
Here $\Gamma_l$'s form a partition of $\Gamma$.

 \begin{rem} With integration by parts, the formulation can accommodate  bases for overlapping domains; e.g.,  piecewise linear polynomials   in \cite{khodayi2019varnet}, which % are used as in the classical finite element methods.  The piecewise linear polynomials 
form a complete basis in $H^1$. 
\end{rem}

\subsection{Error estimates}
In this section, we present error estimates for continuous hp-VRM but do not present detailed analysis for the discrete \textit{hp}-VRM, as it is straightforward to combine
the analysis of 
discrete RM and continuous \textit{hp}-VRM.

% As derived in Section \ref{sec:pinn-error}, 
The key idea   is to identify 
a class of functions and a set of orthogonal basis
that satisfy a certain norm relation between 
the projection operator
and the full operator.
We then apply Theorem~\ref{thm:convergence-nn-linear-prior-posteriori-generic}
to derive an error estimate 
for the hp-VRM.

For a set $\{\Phi_{k,i}\}_{i=1}^{N_k}$  of orthogonal basis for $Y|_{\Omega_k}$, $k=1,\cdots, K$,
let us define the associated projection operator $P_{h,\bm{N}}$ by
\begin{equation} \label{def:projection}
    P_{h,\bm{N}}v= \displaystyle\sum_{k=1}^{K} 
    P_{k} v, \quad \text{where} \quad 
    P_{k}v = \sum_{i=1}^{N_k} (v, \Phi_{k,i})_{Y|_{\Omega_k}} \Phi_{k,i}, \qquad \forall v \in Y,
\end{equation}
with $\bm{N} = (N_1,\cdots, N_K)$.

Following the idea in Section
\ref{sec:ssec-effect-discretizing-integrals}, we may
assume that  
there exists a compact set $Y_c$ of $Y$ such that 
for all $v\in \mathfrak{N}_{\theta,n}\cap X\cap Y_c$ such that 
$2\norm{P_{h,\mathbf{N}}v}_Y\geq \norm{P_{h,\mathbf{N}}v}_Y$. As in Section \ref{sec:ssec-effect-discretizing-integrals}, such inequality relies on 
Bernstein-type inequality for networks but is unavailable for deep feed-forward networks.  Instead, we follow 
a similar approach used in Section  \ref{sec:prototype-PINNs}.

%
%\vskip-20pt 
%

\begin{defn}[Definition of $\tilde{V}_{K}$] \label{def:V_tilde_K}
    Let Assumptions~\ref{assum:solution-existence} 
    {and~\ref{assum:approximation-abilty-X}} hold.
    Let $\{u_n^*\}$ be an approximation sequence in $X$ 
    from Definition~\ref{def:solution}.
    {For a positive decreasing sequence $\{\epsilon_n\}$ that converges to 0,
    let $u_{m_n} \in \mathfrak{N}_{\theta,m_n}$ be a neural network 
    satisfying $\|u_{m_n} - u_n^*\|_X \le \epsilon_n$ (this is guaranteed by Assumption~\ref{assum:approximation-abilty-X}).}
    Let $\Omega = \bigcup_{k=1}^K \Omega_k$.
    For each $k$,
    let $G_{k}$ be a compact set in $Y|_{\Omega_k}$
    containing the sequence {$\{ (Au_{m_n} - f)|_{\Omega_k} \}_{n\ge 1}$}.
%     that depends only on $K$ and {\color{red}	$\{u_{m_n}\}$}.
    We then define a class of functions in $X$ 
    as follows:
    \begin{equation*}
        \tilde{V}_{K}:=\{v \in X : (Av - f)|_{\Omega_k} \in G_{k}, \forall k = 1,\cdots,K \}.
    \end{equation*}
\end{defn}

Next, we show that 
the function class $\tilde{V}_K$ 
is sufficiently large enough to reach 
a zero training loss.

%----------------------------------------------
\begin{prop}[\textbf{Loss Convergence}] \label{prop:vpinn-loss-convg}
	Suppose 
	Assumptions~\ref{assum:solution-existence}, 
	~\ref{assum:approximation-abilty-X} 
	and
	{
		\eqref{eq:norm-equivalence-generic-B}
		of Assumption~\ref{assum:norm-equivalence-generic} hold.}
	For any $\bm{N}=(N_1,\dots,N_K)$ and $\tau \ge 1$,
	let $\mathcal{J}_\tau^{h,\bm{N}}$ be the loss functional 
	\eqref{eq:loss-function-hp-vpinn}
	and $u_{\mathbb{N},n}^{\tau,\bm{N}}$ be its quasi-minimizer
	\footnote{It means that 
		$\mathcal{J}_\tau^{h,\bm{N}}(u_{\mathbb{N},n}^{\tau,\bm{N}})
		\le \inf_{v \in \mathfrak{N}_{\theta,n}\cap \tilde{V}_K} \mathcal{J}_\tau^{h,\bm{N}}(v) + \epsilon$, where $\epsilon\geq 0$ is small.}.
%	of $\mathcal{J}_\tau^{h,\bm{N}_{\epsilon}}$ \eqref{eq:minimization-prob-vpinn}. 
	Then,
	$\lim_{n\to \infty} \mathcal{J}_\tau^{h,\bm{N}}(u_{\mathbb{N},n}^{\tau,\bm{N}}) = 0$.
\end{prop}
\begin{proof}
	{ 
		Let $\{u_k^*\}$ and $\{u_{n_k}\}$ be
		the sequences from Definition~\ref{def:V_tilde_K}. 
		Let $u_{\mathbb{N},n}^{\tau,\bm{N}} \in \mathfrak{N}_{\theta,{n}}\cap \tilde{V}_K$ be a quasi-minimizer of the loss functional,
		i.e.,
		$ % \begin{equation*}
			\mathcal{J}_\tau^{h,\bm{N}}(u_{\mathbb{N},n}^{\tau,\bm{N}})
			\le \inf_{v \in \mathfrak{N}_{\theta,n}\cap \tilde{V}_K} \mathcal{J}_\tau^{h,\bm{N}}(v) + \delta_n,
			$ %  \end{equation*}
		where $\lim\limits_{n \to \infty}\delta_n= 0$.
		Since $u_{n_k} \in \mathfrak{N}_{\theta,{n_k}}\cap \tilde{V}_K$,
		we have
		\begin{align*}
			\mathcal{J}_\tau^{h,\bm{N}}(u_{\mathbb{N},n_k}^{\tau,\bm{N}})
			&\le 
			\mathcal{J}_\tau^{h,\bm{N}}(u_{n_k}) + \delta_{n_k}
			\le
			\mathcal{J}_\tau(u_{n_k}) + \delta_{n_k}.
		\end{align*}
		Since $\lim_{k\to \infty} \mathcal{J}_\tau(u_{n_k}) = 0$
		(by \eqref{eq:norm-equivalence-generic-B}
		of Assumption~\ref{assum:norm-equivalence-generic}),
		we have 
		$\lim_{k\to \infty} \mathcal{J}_\tau^{h,\bm{N}}(u_{\mathbb{N},n_k}^{\tau,\bm{N}}) = 0$.
		Since
		$\mathcal{J}_\tau^{h,\bm{N}}(u_{\mathbb{N},s_1}^{\tau,\bm{N}}) \le \mathcal{J}_\tau^{h,\bm{N}}(u_{s_2}) + \delta_{s_1}$
		for all $s_1 < s_2$,
		the proof is completed.}
\end{proof}

%Proposition~\ref{prop:vpinn-loss-convg}
%indicates that the function class $\tilde{V}_K$
%is sufficiently large enough to reach 
%a zero training loss.
The training loss being zero, however, does not necessarily 
imply the convergence of quasi-minimizers to 
the solution to Equation  \eqref{eq:linear-generic-bvp}.
%In the next lemma, we show how the hp-VRM loss is 
%related to the approximation error.
In the next lemma, we show that if the hp-VRM loss is carefully constructed, the loss functional of \eqref{eq:loss-function-hp-vpinn}
is a good approximation to the loss functional of
\eqref{eq:loss-functional}.

% To establish the norm relations, we need to impose some requirements on the projection to be used. 
{
\begin{lem} \label{lem:projection-full}
    Let Assumptions~\ref{assum:solution-existence} 
    {and~\ref{assum:approximation-abilty-X}} hold.
    Let $\tilde{V}_{K} \subset X$
    be a class of functions defined in 
    Definition~\ref{def:V_tilde_K}.
    Let $\Omega = \bigcup_{k=1}^K \Omega_k$.
%    $\{u_n^*\}$ be a solution sequence in $X$ 
%    from Definition~\ref{def:solution},
%    and
    For any $\epsilon > 0$,
    there exists a set of orthogonal basis 
    $\{\Phi_{k,i}\}_{i=1}^{N_{\epsilon,k}}$
    % $\{\Phi_{k,1},\dots,\Phi_{k,N_{\epsilon,k}}\}$
    in $Y|_{\Omega_k}$ for $k=1,\dots,K$
    that defines $\mathcal{J}_{\tau}^{h,\bm{N}_{\epsilon}}$
    with  $\bm{N}_{\epsilon}=(N_{\epsilon,1},\cdots,N_{\epsilon,K})$,
    which satisfies 
    %a set $\{\Phi_{k,i}\}_{i=1}^{N_{\epsilon,k}}$ of orthonormal basis with respect to $Y|_{\Omega_k}$
    %that depends on $\epsilon$ and $K$ such that 
  %  \begin{align*}
  $  \mathcal{J}_{\tau}^{h,\bm{N}_{\epsilon}} (v)
    \ge \mathcal{J}_{\tau} (v) - \epsilon, \forall v \in \tilde{V}_{K}.
   $ %  \end{align*}
\end{lem}
\begin{proof}
%    Since $\{u_n^*\}_{n \ge 1}$ is a sequence 
%    from Definition~\ref{def:solution},
%    the set $\{Au_n^* - f\}_{n \ge 1}$
%    is relatively compact in $Y$.
    Let $f_k = f|_{\Omega_k}$,
    and let $G_k$ be the compact set in $Y|_{\Omega_k}$
    from Definition~\ref{def:V_tilde_K}.
%    containing $\{Au_{m_n}|_{\Omega_k} -f_k\}_{n \ge 1}$.
%    Note that there exists
%    at least one such compact set,
%    namely, $\{Au_{n}^*|_{\Omega_k} -f_k\}_{n \ge 1} \cup \{0\}$.
    Then, for any $\epsilon > 0$, 
    there exists a finite-dimensional subspace $\tilde{K}_{Y,k}^{\epsilon,K}$ of $Y|_{\Omega_k}$
    such that for any $u \in G_{k}$,
    there exists $u_{\epsilon} \in \tilde{K}_{Y,k}^{\epsilon,K}$
    satisfying $\|u_{\epsilon} - u\|_{Y|_{\Omega_k}} \le \sqrt{\epsilon/K}$.
    Let $\tilde{K}_{Y,k}^{\epsilon,K}$ be spanned by 
    $\{\Phi_{k,i}\}_{i=1}^{N_{\epsilon, k}}$
    and let $P_{h,\bm{N}_{\epsilon}}$ be defined through 
    this basis.
%     $\{\Phi_{k,i} : i=1,\dots, N_{\epsilon, k}\}_{k=1}^K$.
	Observe that for any $v \in \tilde{V}_K$,
	since $Av|_{\Omega_k} - f_k \in G_k$,
	there exists $\tilde{g}_k \in \tilde{K}_{Y,k}^{\epsilon,K}$
	such that $\|Av|_{\Omega_k} - f_k - \tilde{g}_k \|_{Y|_{\Omega_k}} \le \sqrt{\epsilon/K}$.
	Hence, $\|(I-P_{N_{\epsilon,k}})(Av|_{\Omega_k} - f_k)\|_{Y|_{\Omega_k}} \le \sqrt{\epsilon/K}$.
    It then can be checked that 
%    for any $v \in \tilde{V}_K$,
    $
        \|(I - P_{h,\bm{N}_{\epsilon}})(Av - f)\|_Y^2
        = \sum_{k=1}^K 
        \|(I - P_{N_{\epsilon,k}})(Av|_{\Omega_k} - f_k)\|_{Y|_{\Omega_k}}^2
        \le \epsilon, \forall v \in \tilde{V}_k.
    $ %\end{align*}
    Thus, we have 
     $%\begin{align*}
    \mathcal{J}_{\tau}^{h,\bm{N}_{\epsilon}} (v)
    =
    \mathcal{J}_{\tau} (v)
    - \|(I - P_{h,\bm{N}_{\epsilon}})(Av - f)\|_Y^2
    \ge \mathcal{J}_{\tau} (v) - \epsilon,
  $ %  \end{align*}
    % By letting $G_k = \tilde{K}_{Y,k}^{\epsilon,K}$,
    and the proof is completed.
\end{proof}
}

To establish error estimates and convergence,
we now consider 
the hp-VRM formulation under 
$\mathfrak{N}_{\theta,n} \cap \tilde{V}_{K}$
{with the loss functional defined through the specific basis from Lemma~\ref{lem:projection-full}.}
That is,
\begin{equation} \label{eq:minimization-prob-vpinn}
	\min_{v \in \mathfrak{N}_{\theta,n} \cap \tilde{V}_{K} } \mathcal{J}_{\tau}^{h,\bm{N}_\epsilon}(v).
\end{equation}
The feasible function class $\tilde{V}_K$
depends only on the number of partitions of the domain $\Omega$
and an approximation sequence {$\{u_{m_k}\}_{k \ge 1}$ from Definition~\ref{def:V_tilde_K}.}
Also, $\tilde{V}_K$ is non-empty as it includes $\{u_{m_k}\}_{k \ge 1}$.

%----------------------------------------------
\begin{thm}[\textbf{Error estimates for hp-VRM}] 
\label{thm:prior-posteriori-generic-hp-VPINN}
Under the same assumptions in Lemma~\ref{lem:projection-full},
suppose Assumption~\ref{assum:norm-equivalence-generic} holds.
Let $\tilde{V}_{K} \subset X$ be a function class defined in Definition~\ref{def:V_tilde_K}.
For any $\epsilon > 0$ and $\tau \ge 1$,  
let $\{\Phi_{k,i}\}_{i=1}^{N_{\epsilon,k}}$ be a set of orthonormal basis with respect to $Y|_{\Omega_k}$ for $k=1,\dots, K$ from Lemma~\ref{lem:projection-full}.
% and that defines the loss functional $\mathcal{J}_{\tau}^{h,\bm{N}_{\epsilon}}$ \eqref{eq:loss-function-hp-vpinn}. 
Let $u_{\mathbb{N},n}^{\tau,\bm{N}_{\epsilon}} \in \mathfrak{N}_{\theta,n} \cap \tilde{V}_{K}$ be a quasi-minimizer
of \eqref{eq:minimization-prob-vpinn},
and $u^*$ be the solution to \eqref{eq:linear-generic-bvp}
in the sense of Definition~\ref{def:solution}.
Then, the following a posterior estimation holds:
% for any $\epsilon' > 0$, there exists $u_{\epsilon'}^* \in X$ such that
\vskip -18pt 
\begin{eqnarray*} 
    \|u_{\mathbb{N},n}^{\tau,\bm{N}_{\epsilon}} - u^*\|_V
    \le 
	\sqrt{2}C_1^{-1}\left( \mathcal{J}_{\tau}^{h,\bm{N}_{\epsilon}}(u_{\mathbb{N},n}^{\tau,\bm{N}_{\epsilon}}) + \delta_n +
    \epsilon
    \right)^{\frac{1}{2}}, 
\end{eqnarray*}
% \vskip -5pt  \noindent
where $C_1$ is the   constant defined in the norm relation 
\eqref{eq:norm-equivalence-generic}
and $\delta_n$ is a vanishing sequence stemming from 
the choice of quasi-minimizers.
\end{thm}
\begin{proof}
Let $u^*$ be the solution to \eqref{eq:linear-generic-bvp}.
%Let $\{u_n^*\}$ and $\{u_{m_n}\}$ be
%the sequences from Definition~\ref{def:V_tilde_K}.
%and $\{u_n^*\}_{n \ge 1}$ be its corresponding sequence in $X$
%from Definition~\ref{def:solution}
%(also Assumption~\ref{assum:solution-existence}).
For $\epsilon > 0$, let $u_{\mathbb{N},n}^{\tau,\bm{N}_{\epsilon}}$ be 
a quasi-minimizer of the loss $\mathcal{J}_{\tau}^{h,\bm{N}_{\epsilon}}$ 
defined through Lemma~\ref{lem:projection-full} over $\mathfrak{N}_{\theta,n} \cap \tilde{V}_{K}$. 
That is,
$\mathcal{J}_{\tau}^{h,\bm{N}_{\epsilon}}(u_{\mathbb{N},n}^{\tau,\bm{N}_{\epsilon}}) \le \inf_{w\in \mathfrak{N}_{\theta,n}\cap \tilde{V}_{K}}   \mathcal{J}_{\tau}^{h,\bm{N}_{\epsilon}}(w) + \delta_n$,
%\begin{align*}
%    \mathcal{J}_{\tau}^{h,\bm{N}_{\epsilon, K}}(u_{n}^{\tau,\bm{N}_{\epsilon, K}}) \le \inf_{w\in \mathfrak{N}_{\theta,n}\cap \tilde{V}_{\epsilon, K}}   \mathcal{J}_{\tau}^{h,\bm{N}_{\epsilon, K}}(w) + \delta_n,
%    % \le 
%    % \mathcal{J}_{\tau}^{h,\bm{N}}(u_{k_n}^*) + \delta_n,
%\end{align*}
%\vskip -11pt \noindent
where  $\delta_n \to 0$ as $n \to \infty$.
By Lemma~\ref{lem:projection-full}, 
since $\mathcal{J}_{\tau}(v) \le \mathcal{J}_{\tau}^{h,\bm{N}_{\epsilon}}(v) + \epsilon$ 
for all $v \in \tilde{V}_{K}$, 
we obtain
\begin{align*}
    \mathcal{J}_{\tau}(u_{\mathbb{N},n}^{\tau,\bm{N}_{\epsilon}})
    &= \mathcal{J}_{\tau}^{h,\bm{N}_{\epsilon}}(u_{\mathbb{N},n}^{\tau,\bm{N}_{\epsilon}}) + \|(I-P_{h,\bm{N}_{\epsilon}})(Au_n^\tau -f)\|_Y^2  
    \le \inf_{w\in \mathfrak{N}_{\theta,n}\cap \tilde{V}_{K}}  \mathcal{J}_{\tau}^{h,\bm{N}_{\epsilon}}(w) + \delta_n +
    \epsilon.
    % 2(\epsilon + \|I-P_{h,\bm{N}_{\epsilon, K}}\|^2_{op} \inf_k \|Au_n^\tau - Au_{k}^*\|_Y^2 ).
\end{align*}
It then follows from Theorem~\ref{thm:convergence-nn-linear-prior-posteriori-generic} that 
\begin{equation*}
    \begin{split}
    C_1^2  \|u_{\mathbb{N},n}^{\tau,\bm{N}_{\epsilon}} - u^*\|_V^2  \le  
    2 \mathcal{J}_{\tau}(u_{\mathbb{N},n}^{\tau,\bm{N}_{\epsilon}})  
    \le 
	2\left(\inf_{w\in \mathfrak{N}_{\theta,n}\cap \tilde{V}_{K}}  \mathcal{J}_{\tau}^{h,\bm{N}_{\epsilon}}(w) + \delta_n +
    \epsilon
    \right).
    \end{split}
\end{equation*}
% Note that 
% \begin{align*}
%     \inf_k \|Au_n^{\tau,\bm{N}} - Au_k^*\|_Y^2
%     \le 2\|Au_n^{\tau,\bm{N}} - Au^*\|_Y^2
%     + 2\inf_k \|Au_k^* - Au^*\|_Y^2
%     = 2\|Au_n^{\tau,\bm{N}} - f\|_Y^2.
% \end{align*}
Letting $w = u_{\mathbb{N},n}^{\tau,\bm{N}_{\epsilon}}$ completes the proof. 
\end{proof}

\begin{thm}
    [Convergence of hp-VRM] \label{thm:convg-hpVPINNs}
    Under the same conditions and assumptions of Theorem~\ref{thm:prior-posteriori-generic-hp-VPINN}, 
    for $\epsilon > 0$, let $u_{\mathbb{N},n}^{\tau,\bm{N}_{\epsilon}} \in \mathfrak{N}_{\theta,n} \cap \tilde{V}_{K}$ be a quasi-minimizer of \eqref{eq:minimization-prob-vpinn}
    defined through Lemma~\ref{lem:projection-full}. Then,  $ %  \begin{equation*} 
    \lim_{\epsilon \to 0} \lim_{n\to \infty} 
    \|u_{\mathbb{N},n}^{\tau,\bm{N}_{\epsilon}} - u^*\|_V
    = 0, 
    $ %  \end{equation*}
    where $u^*$ is the solution to \eqref{eq:linear-generic-bvp}
    from Assumption~\ref{assum:solution-existence}. 
 
\end{thm}
\begin{proof}
    It follows from
    Theorem~\ref{thm:prior-posteriori-generic-hp-VPINN}
    and Proposition~\ref{prop:vpinn-loss-convg}
    that
$%    \begin{equation*}     \begin{split}
    \lim_{n\to \infty} 
    \|u_{\mathbb{N},n}^{\tau,\bm{N}_{\epsilon, K}} - u^*\|_V
    \le \sqrt{2}C_1^{-1}
    \epsilon^{\frac{1}{2}}.
$ %    \end{split}     \end{equation*}
    By letting $\epsilon \to 0$, 
    the proof is completed.
\end{proof}

\begin{rem}[Weaker formulation]
Let $\set{\Phi_{k,i}}$ be an orthonormal basis in $H^{1}_0(\Omega)$. Then the functional in the hp-VRM formulation \eqref{eq:loss-function-hp-vpinn} becomes
\begin{equation}
\mathcal{J}_\tau  (v) =  \norm{f- A v}^2_{H^{-1}}+  \tau \norm{Bu-g}^2_{Z},\quad \tau\geq 0,
\end{equation} 
when $h\to 0$ and $N\to\infty$.
If  $\nabla \Phi_{k,i}$   form a  complete basis in $L^2$ (which is the case for piecewise linear polynomials), 
then  performing integration-by-parts  will lead to problems as in the case of $L^2$ orthonormal bases. For example,  when $A v=-\Delta v$, we then have 
\begin{equation}
\mathcal{J}_\tau  (v) =  \norm{F- \nabla v}^2_{L^2}+  \tau \norm{Bu-g}^2_{Z},\quad \tau\geq 0.
\end{equation} 
Here $f= {\rm div } \, F$ and $F\in [H^1_0(\Omega)]^2$ and we assume $\nabla\Phi_{k,i}$ form a complete orthogonal basis. In fact, 
$\sum_{k,i=1}^\infty(f-\Delta v,\Phi_{k,i})=\sum_{k,i}^\infty (-F + 
\nabla v, \nabla \Phi_{k,i}) =\norm{\nabla v- F}^2_{L^2}$. 
\end{rem}

%% file: discussion.tex
\section{Conclusion and Discussions}
We proposed an abstract framework for analyzing the convergence of residual minimization for linear PDEs using neural networks. 
When Bernstein-type inequalities are available for neural networks, we use the discrete norm relations to obtain the convergence; See Theorem \ref{thm:convergence-nn-linear-prior-posteriori-generic-discrete} and Example \ref{exm:bernstein-ineq-guassian-discrete-norm-relation}. 
When Bernstein-type inequalities are unavailable, we use the  Rademacher complexity to obtain the convergence. We also present some examples in \ref{sec:elliptic-advection-fractional} on verification of our assumptions. 
Both approaches introduced tailored classes of neural networks
that enjoy some desired properties.
% as in nonlinear problems.

The framework developed in this paper may serve as guidance in designing loss functionals of the residual minimization. First, we need the stability of the equations under user-defined metrics. Second, we need to balance the number of training points and the networks' size as in Example \ref{exm:bernstein-ineq-guassian-discrete-norm-relation}.
%
%
%
%Though the abstract framework can explain the convergence of networks for some problems, such as Poisson equations,
%it is not readily applicable for problems with nonlinearities.
However, the conditions for convergence are 
  not readily verifiable for deep neural networks.
Besides the limitations mentioned in Section \ref{sec:intro}, 
the verification of  Bernstein-type inequalities and the Rademacher complexity is limited to two-layer networks. Also, the verification may be complicated 
depending on the operators and equations under consideration.  
These aspects are being investigated in the community, and more efforts are required.  

%% file: appendix.tex
\appendix

\section{Proof of Theorem 4.3} \label{APP:THM:INTP-DERIV-NNs}
The proof is similar to the proof 
  of Proposition 4.5 in \cite{Mhaskar06Markov}. The key step is to establish the 
  following lemma.

  \begin{lem}\label{eq:interpolate-gaussian-rbf-estimate}Let $g(x)= \exp \left(-\abs{x- {w}}^{2} \right)  $ where $x,w\in\Real^d$. Let
  	$\mathcal{I}_n$ be the interpolation operator defined in \eqref{eq:weighted-interpolation}. Then there exist constant $c_1,c>0$ such that
  	\[	 \norm{\mathbf{D}^{\mathbf{j}}g  -  \mathcal{I}_n\mathbf{D}^{\mathbf{j}}g}_{L^2}\leq  c_1 n^c  
  	\frac{(\sqrt{2 d}\abs{ {w}})^{n+1} \exp \left(d\abs{w}^{2}\right)}{\sqrt{n !}},\quad \text{for any multiindex } \mathbf{j}\in \mathbb{N}^d.\]	
  \end{lem}
  
  The observation here is that $g(x)$ is closely related to  Hermite polynomials. 
  Let $h_k(x_1)$ be the orthonormal Hermite polynomials on the real line with respect to the weight $\exp(-\abs{x_1}^2))$: $\int_{\Real} h_k(x_1)h_j(x_1)\exp(-\abs{x_1}^2)\,dx_1=\delta_{k,j}$.  
  Then by the generating function of the Hermite polynomials, we have 
  $\exp(2x_1t-t^2) = \pi^{1/4}\sum_{k=0}^\infty \frac{h_k(x_1)}{\sqrt{k!}}(\sqrt{2}t)^k$
  For the multiindex $\mathbf{k}$, $h_{\bf k} =\Pi_{j=1}^s h_{k_j}(x_j)$. 
  Then with the standard multivariate notation, we have 
  $$g(x,w)= \exp(-\abs{x-w}^2) =\pi^{d/4}\sum_{\abs{\mathbf{k}}\geq 0} \frac{h_{\mathbf k}(x)\exp \left(-\abs{x}^{2}\right)}{\sqrt{\mathbf{k}!}}(\sqrt{2}w)^{\mathbf{k}}. $$

  For integer $n \geq 1$, let
  \begin{equation}\label{eq:truncated-generating-function-weighted}
  P_{n}( {x},  {w}):=\pi^{d / 4} \sum_{0 \leq|\mathbf{k}| \leq n} \frac{h_{\mathbf{k}}( {x}) \exp \left(-\abs{x}^{2}\right)}{\sqrt{\mathbf{k} !}}(\sqrt{2}  {w})^{\mathbf{k}}, \quad 
  P_{n}^\perp( {x},  {w}) =g( {x},  {w})- P_{n}( {x},  {w}).
  \end{equation} 
  Then   by Lemma 4.6 of \cite{Mhaskar06Markov},   there exist constants $c_1,c>0$ such that
  \begin{equation}\label{eq:approximate-gaussian-rbf-by-hermite}
   	\left\|\exp \left(-\abs{\circ-{w}}^{2}\right)-P_{n}(\circ,  {w})\right\|_{W^{r,2} } \leq c_{1} n^{c} \frac{(\sqrt{2 d}\abs{ {w}})^{n+1} \exp \left(d\abs{w}^{2}\right)}{\sqrt{n !}},\quad \text{for any integer } r,n\geq 1.
  \end{equation}

  \begin{proof}
  	By the fact that $I_N^c (\exp(\frac{\abs{x}^2}{2})P_n(x,w)) = \exp(\frac{\abs{x}^2}{2})P_n(x,w)$ ($N\geq n$) and the triangle inequality, 
  	\begin{eqnarray*}
  		&&\norm{\mathbf{D}^{\mathbf{j}}g  -  \exp(-\frac{\abs{\circ}^2}{2}) I_N^c (\exp(\frac{\abs{\circ}^2}{2})\mathbf{D}^{\mathbf{j}}g  }_{L^2}\\
  		&\leq&
  		\norm{\exp(\frac{\abs{\circ}^2}{2})\mathbf{D}^{\mathbf{j}}g   -I_N^c \left(\exp(\frac{\abs{\circ}^2}{2}\mathbf{D}^{\mathbf{j}}g ) \right)}_{L^2} 
  		\\
  		&\leq &  	\norm{\exp(\frac{\abs{\circ}^2}{2}) \mathbf{D}^{\mathbf{j}}(g -P_n(\circ,w) )   }_{L^2} +  \norm{I_N^c \left(\mathbf{D}^{\mathbf{j}}(g -P_n(\circ,w) )\exp(\frac{\abs{\circ}^2}{2})\right)}_{L^2}=: I + II.
  	\end{eqnarray*}

  We first estimate the term $I$.  	
  	A careful check of the proof of Lemma 4.6 in \cite{Mhaskar06Markov} leads to 
  	\begin{equation}\label{eq:approximate-gaussian-rbf-by-hermite-v1}
  	   I=\norm{\exp(\frac{\abs{\circ}^2}{2})\mathbf{D}^{\mathbf{j}} P_n^\perp(\circ,w) }_{L^2}\leq  c_1 n^c
  	\frac{(\sqrt{2 d}\abs{ {w}})^{n+1} \exp \left(d\abs{w}^{2}\right)}{\sqrt{n !}}.
  	\end{equation}
  We now estimate the term $II$.  By the stability of the interpolation operator 
  $I_N^c$ (multi-dimensional analogue of Lemma 3.1 in  \cite{GuoSX03}),
  \begin{eqnarray*}
   \norm{I_N^c \left(\mathbf{D}^{\mathbf{j}} P_n^\perp(\circ,w)  \exp(\frac{\abs{\circ}^2}{2})\right)}_{L^2}   &\leq& c\sum_{\abs{\mathbf{i}}\leq 1} \norm{\mathbf{D}^{\mathbf{i}}\left(\mathbf{D}^{\mathbf{j}}P_n^\perp (\circ,w)   \exp(\frac{\abs{\circ}^2}{2})\right) }_{L^2}.
  \end{eqnarray*}	
%Then by \eqref{eq:truncated-generating-function-weighted} and the generation function representation of $g$, we have 
 %\begin{equation}\label{eq:approximate-gaussian-rbf-residual}
 %  g -P_n(x,w) = \pi^{d / 4} \sum_{\abs{\mathbf{k}} \geq n+1} \frac{h_{\mathbf{k}}( {x}) \exp \left(-\abs{x}^{2}\right)}{\sqrt{\mathbf{k} !}}(\sqrt{2}  {w})^{\mathbf{k}}.
 %\end{equation} 
 %
 Recall that in Theorem 2.2 of \cite{GuoSX03} \label{eq:inverse-inequality}, it is shown that for any $v =Q_N \exp(-\frac{\abs{x}^2}{2})$, $x\in\Real$, then for $m\geq 0$, $\norm{\partial_x^m v}_{L^2}\leq  (2N+1)^{\frac{m}{2}}\norm{v}_{L^2}$. Since we are using a tensor product of interpolation, a multi-dimensional inverse estimate also holds $\norm{\mathbf{D}^{\mathbf{i}}\left(h_{\mathbf{k}}  \exp (-\abs{\circ}^{2}/2)\right)}_{L^2}\leq C (\abs{\mathbf{k}}+1)$ for $\abs{\mathbf{i}}\leq 1$. Then, by \eqref{eq:truncated-generating-function-weighted} and the inverse estimate,	\begin{eqnarray*}
	&&\norm{I_N^c \left(\mathbf{D}^{\mathbf{j}}  P_n^\perp (\circ,w) \exp(\frac{\abs{x}^2}{2})\right)}_{L^2} \\	&\leq& 
	c\sum_{\abs{\mathbf{i}}\leq 1}  \norm{\mathbf{D}^{\mathbf{i}}\Big(\mathbf{D}^{\mathbf{j}} P_n^\perp (\circ,w) \exp(\frac{\abs{\circ}^2}{2})\Big) }_{L^2}  \\
	&\leq& c\sum_{\abs{\mathbf{i}}\leq 1}  \sum_{|\mathbf{k}| \geq n+1} \norm{\mathbf{D}^{\mathbf{i}} \left( h_{\mathbf{k}+\mathbf{j}}(\circ) \exp  (-\frac{\abs{\circ}^{2}}{2})\right)}_{L^2}
	\abs{(-\sqrt{2})^{|\mathbf{j}|} \frac{\sqrt{(\mathbf{k}+\mathbf{j}) !}}{\mathbf{k} !} 
		(\sqrt{2} \mathbf{w})^{\mathbf{k}} }\\
	&\leq &   c\sum_{\abs{\mathbf{i}}\leq1}  \sum_{|\mathbf{k}| \geq n+1}( \abs{\mathbf{k}+\mathbf{j}}+1)
	\sqrt{2})^{|\mathbf{j}|} \frac{\sqrt{(\mathbf{k}+\mathbf{j}) !}}{\mathbf{k} !} 
		\abs{\sqrt{2} \mathbf{w})^{\mathbf{k}} }.%\\
%	&\leq& c \sum_{|\mathbf{k}| \geq n+1}  (\abs{\mathbf{k}}+\abs{\mathbf{j}})^{d/2}
%	\abs{(-\sqrt{2})^{|\mathbf{j}|} \frac{\sqrt{(\mathbf{k}+\mathbf{j}) !}}{\mathbf{k} !} 
%		(\sqrt{2} \mathbf{w})^{\mathbf{k}} }\\
\end{eqnarray*}	
In the second inequality, we use the fact  (by  Rodrigues’ formula) that 
\begin{equation*}
  \mathbf{D}^{\mathbf{j}} \left(h_{\mathbf{k}} \exp(-\abs{x}^2)\right) = (-\sqrt{2})^{\mathbf{j}}  \frac{\sqrt{(\mathbf{k}+\mathbf{j}) !}}{
  \sqrt{\mathbf{k} !}} h_{\mathbf{k}+\mathbf{j}} \exp(-\abs{x}^2)
\end{equation*}
 Then by the same argument in the proof of Lemma 4.6 in \cite{Mhaskar06Markov},  there exist $c_1,c>0$ such that 
  \begin{equation}\label{eq:stability-interpolation}
  II=\norm{I_N^c \left(\mathbf{D}^{\mathbf{j}} P_n^\perp(\circ,w)  \exp(\frac{\abs{\circ}^2}{2})\right)}_{L^2} \leq  c_1  n^c
  \frac{(\sqrt{2 d}\abs{ {w}})^{n+1} \exp \left(d\abs{w}^{2}\right)}{\sqrt{n !}}.
  \end{equation}
  Then the desired conclusion follows from \eqref{eq:stability-interpolation} and 
 \eqref{eq:approximate-gaussian-rbf-by-hermite-v1}.  
  \end{proof}

  The proof of Theorem \ref{THM:INTP-DERIV-NNs} is verbatim the same as that  of Proposition 4.5  in \cite{Mhaskar06Markov}, except using Lemma 
 \ref{eq:interpolate-gaussian-rbf-estimate} in place of  the estimate \eqref{eq:approximate-gaussian-rbf-by-hermite}.

%% file: elliptic.tex
%%%%%%%%%%%%%%%%%%%%%%%%%%%%%%%%%%%%%%%%%%%%%%%%%  
\section{Illustration Examples and Verification of Assumptions}	 \label{sec:elliptic-advection-fractional}
%%%%%%%%%%%%%%%%%%%%%%%%%%%%%%%%%%%%%%%%%%%%%%%%%

In this section, 
we consider linear elliptic, advection equations, and an integro-differential equation and demonstrate the  key assumptions in Section \ref{sec:pinn-error} are satisfied.

 \subsection{Elliptic problems} \label{sec:elliptic}
Let $A$ be a linear differential operator of the form
 \begin{equation} 
 	 A= -\sum_{i,j=1}^d a_{i,j}(x) \partial_{x_i}\partial_{x_j} + \sum_{i=1}^d b_i(x)\partial_{x_i} + c(x),
 	 \qquad
 	 \partial_{x_j} := \frac{\partial}{\partial x_j},
 \end{equation}
and $B$ be the identity operator, i.e., $B={\rm Id}$, which leads to a Dirichlet boundary condition on $\partial\Omega$. 
Also, we make the following assumptions.
 \begin{assu}\label{assu:elliptic-operator}  
 The coefficients $a_{i,j}$ satisfy the uniformly elliptic condition and the   coefficients are   in $C^2$, i.e., twice continuously differentiable.
 Also, the only solution with zero input data is the zero solution.
 \end{assu}

\begin{lem}[Theorem 2.1 in \cite{BraSch70}] \label{lem:norm-relation-proof-elliptic} 
In addition to Assumption \ref{assu:elliptic-operator}, assume that $A$ is     with $C^\infty$ coefficients, defined on $C^{\infty}$ bounded domain $\Omega$ on $\Real^d$. For any real number $l$, 
 $ %	 \begin{equation}
 	 \norm{u}_{H^l(\Omega)} \leq  C (\norm{Au}_{H^{l-2}(\Omega)}+\abs{u}_{H^{l-1/2}(\partial\Omega)} ), 
 $ % 	 \end{equation} 
 	 for all $u\in C^\infty(\bar{\Omega})$ and $C$ is independent of $u$. 
\end{lem}
 From this lemma, we can deduce from the density argument that 
 for all $l\leq 1/2$, 	 \begin{equation}
\norm{u}_{H^l(\Omega)} \leq  C (\norm{Au}_{H^{-3/2}(\Omega)}+\norm{u}_{L^2(\partial\Omega)} )\leq C (\norm{Au}_{L^2(\Omega)}+\norm{Bu}_{L^2(\partial \Omega)} ),
\end{equation} 
holds for any $u\in H^2(\Omega)$.
% \textcolor{red}{
% such that $Au \in H^{-3/2}(\Omega)$.} 
Thus,  
\eqref{eq:norm-equivalence-generic-A}
of Assumption \ref{assum:norm-equivalence-generic}
is verified  with
$V= H^{l}(\Omega)$ for any $l\leq 1/2$, $Y=L^2(\Omega)$, and $Z= L^2(\partial \Omega)$.
% \begin{center}
% % $X = H^2(\Omega)\subseteq V= H^{l}(\Omega)$,
% $V= H^{l}(\Omega)$, $Y=L^2(\Omega)$, and $Z= L^2(\partial \Omega)$.
% \end{center}
It follows from the trace inequality (e.g. Theorem 1.6.6 of \cite{brenner2007mathematical})
and $u \in H^2(\Omega)$
that \eqref{eq:norm-equivalence-generic-B}
of Assumption \ref{assum:norm-equivalence-generic}
is verified with
$X = H^2(\Omega) \subset V$.

The verification of Assumption \ref{assum:approximation-abilty-X} 
is straightforward as 
feed-forward neural networks are universal approximators 
in Sobolev-Hilbert spaces, see e.g., \cite{Mha96}. 
Therefore, Theorem \ref{thm:convergence-nn-linear-prior-posteriori-generic}
provides 
the error estimates in $X= H^2(\Omega)$
under the conditions in Lemma
~\ref{lem:norm-relation-proof-elliptic}. 
\vskip 10pt 

%%%%%%%%%%%%%%%%%%%%%%%%%%%%%

%---------------------------
\iffalse 
In Case II, we can use the norms in H\"older spaces. %
%For example, if 
%$f\in C^{0,\alpha}$, $\alpha\in (0,1)$, we can obgtain that th 
 In fact, we have with the assumptions  
 \fi 
%%%%%%%%%%%%%%%%%%%%%%%%%%%%%%%%%%%%%%5
\begin{rem}[{Non-smooth data}, \cite{BraSch70}]	 \label{rem:non-smooth-data}
If $f\in H^{s}(\Omega)$, $-2\leq s\leq 0$,
is a non-smooth data,  
we may use a mollifier 
$\Phi_h$ such that for some $C>0$ independent of $h$ and $f$, the following holds:
$ %\begin{equation*}
	\norm{\Phi_hf -f }\leq C h^s \norm{f}_{H^{-s}}, \text{ and } \norm{\Phi_h f -f }_{-2}\leq C h^{2+s}\norm{f}_{H^{-s}}.
$ %\end{equation*}

%	For $f\in L^1(\Omega)$, we can simply use the averages $\frac{1}{\abs{\Omega_k}}\int_{\Omega_k} f\,dx$ to approximate $f$ in a small domain $\Omega_k$.  
\end{rem}

%% file: advection.tex
%----------------------------------------------
%
%------------------------------------------------
 
\subsection{Advection-reaction problems}\label{sec:advection}

%%----------
Let $\Omega\subset \Real^d$ be an open bounded domain, with Lipschitz
boundary $\partial\Omega $ oriented by a unit outward normal vector $\mathbf{n}$. 
We consider an advection-reaction problem.
Let $\mathbf{b}$ be a smooth vector field in $\mathbb{R}^d$ such that 
$\mathbf{b} \in L^\infty(\Omega)^d$
and $\nabla \cdot \mathbf{b}\in L^\infty(\Omega)$.
Let the inflow boundary be 
%\begin{equation*}
$    \partial\Omega_{-}= \set{x\in \partial\Omega|\, \mathbf{b}\cdot \boldsymbol{n} <0}$.
%\end{equation*}
For $c \in L^\infty(\Omega)$,
let $A$ be the differential operator defined by
\begin{equation}\label{eq:advection-reaction-inflow}
 Au =\mathbf{b}(x) \cdot\nabla    u  + c(x) u, \quad x\in \Omega,
\end{equation}
and 
$B$ be the identity operator  on the inflow boundary $\partial\Omega_{-}$. This formulation also covers  time-dependent
advection-reactions; See Remark \ref{rem:time-dependent-ar}.

% $f\in L^p(\Omega)$, $p\in [1,\infty)$, and  

% Let us introduce some spaces for the advection-diffusion problem.  
For  $p\in [1,\infty]$, we   define the   graph space 
%\begin{equation}
$G^{p}_{\mathbf{b}}(\Omega) = \set{v\in L^p(\Omega)|\mathbf{b} \cdot \nabla v\in L^p(\Omega)}$,
%\end{equation} 
which is endowed with the norm
$\norm{v}_{G^p_{\mathbf{b} }} =(\norm{v}^2_{L^p} + \norm{\mathbf{b}\cdot\nabla  v}^2_{L^p})^{1/2} $.  
%An equivalent norm is  $\norm{v}_{G^p_{\mathbf{b} }} =(\norm{v}^2_{L^p} + \norm{ \nabla\cdot (\mathbf{b} v)}^2_{L^p})^{1/2} $. 
We also introduce the following space % defined over the boundary:
\begin{equation*}
 L^p_{\abs{\mathbf{b}\cdot \bf{n}}}(\partial\Omega)=\set{v \text{ is measurable on } \partial \Omega \mid  \int_{\partial\Omega} \abs{\mathbf{b}\cdot \bf{n} } \abs{v}^p  \,dx<\infty}.  
  \end{equation*}
 %  Under mild conditions (Lemma \ref{lem:trace-adv-reac}), the  space $ L^p_{\abs{\mathbf{b}\cdot \bf{n}}}(\partial\Omega)$ is the collection of trace of functions in the space $G_{\mathbf{b}}^p(\Omega)$.

By the trace inequality (Lemma~\ref{lem:trace-adv-reac}) and 
Poincare's  inequality (Theorem \ref{thm:poincare}),  
it can be checked that the following condition holds 
\begin{equation}\label{cond:norm-equivalence-lp-ar}
C_1 \norm{v}_X \leq \norm{A v}_Y +\norm{Bv}_Z\leq C_2 \norm{v}_X,  
\end{equation} 
where $X = G^p_{\mathbf{b}}(\Omega)$, $Y=L^p(\Omega)$, $p\in (1,\infty)$,  $Z= L^p_{\abs{\mathbf{b}\cdot \mathbf{n}}}(\partial\Omega)$ and  $p \in (1,\infty)$. 
This verifies Assumption~\ref{assum:norm-equivalence-generic}.
 
\begin{lem}[Trace inequality, Lemma 2.1 of \cite{Can17}]\label{lem:trace-adv-reac}
Let $
\partial\Omega_{+} =  \set{x\in \partial\Omega|\, \mathbf{b}\cdot \boldsymbol{n} >0}$. 	
	Assume that the inflow and outflow boundaries are well-separated: $\overline{\partial\Omega_+} \cap \overline{\partial\Omega_{-}} =\emptyset$. % , where 	  
	Let $p\in (1,\infty)$ and  $g\in L^p_{\abs{\mathbf{b}\cdot \mathbf{n}}}(\partial\Omega)$,  there exists $v_g\in G^p_{\mathbf{b}}(\Omega)$ such that 
	$v_g  = g$ on $\partial \Omega_{+}\cup \partial\Omega_{-}$.
	In other words, the trace of $v_g \in  G^p_{\mathbf{b}}(\Omega)$ exists in $L^p_{\abs{\mathbf{b}\cdot \mathbf{n}}}(\partial\Omega)$ and 
	$ %\begin{equation}
	\norm{g}_{L^p_{\abs{\mathbf{b}\cdot \mathbf{n}}}(\partial\Omega)} \leq C \norm{v_g}_{G_{\mathbf{b}}^p}.
	$ %\end{equation}
\end{lem}

\iffalse  For all $v\in G^p_{\mathbf{b}}(\Omega)$ and $w\in G^{p^\ast}_{\mathbf{b}}(\Omega)$, the integration by parts formula holds:
\begin{equation}
\int_{\Omega} (\mathbf{b}\cdot \nabla v)w + (\mathbf{b}\cdot \nabla w)v + (\nabla \cdot \mathbf{b})  vw \,dx = \int_{\partial\Omega   }(\mathbf{b}\cdot \mathbf{n}) vw\,dx.
\end{equation}
\fi

%Let $G^p_{\mathbf{b},0\pm } =\set{v\in  G^p_{\mathbf{b}}(\Omega)| v(x)=0, \,x \in \partial\Omega_{\pm} }$.   
%For the linear operator in  \eqref{eq:advection-reaction-inflow},  the following Poincare inequality holds. 
\begin{thm}[Poincare inequality, \cite{Can17}]
	\label{thm:poincare}
	Let  $p\in (1,\infty)$.  
	Assume that there exists an  Lipschitz continuous function $\eta(x)$ and a positive constant $\mu_1$ such that 
	\begin{equation}\label{eq:condition-ar-well-posedness-in-Lp}
	c(x) -\frac{1}{p}\nabla \cdot \mathbf{b}(x) - \frac{1}{p} \mathbf{b}(x)\cdot \nabla \eta(x) \geq \mu_1>0, \text{ a. e. } x\in\Omega.
	\end{equation}
	Then for $v\in \set{v\in  G^p_{\mathbf{b}}(\Omega)| v(x)=0, \,x \in \partial\Omega_{-} }$, \quad 
$%	\begin{eqnarray}
	\norm{v}_{L^p}  \leq  C \norm{c v + \mathbf{b}\cdot\nabla v }_{L^p}.%,\\
	%	\norm{v}_{L^{p^\ast}} &\leq& C \norm{c v -\nabla \cdot ( \mathbf{b} v) }_{L^{p^\ast}},\quad  v\in G^{p^\ast}_{\mathbf{b},0+}(\Omega),
$ %	\end{eqnarray}
\end{thm}  

The condition \eqref{eq:condition-ar-well-posedness-in-Lp} 
in Theorem~\ref{thm:poincare}
can be satisfied in the   cases of Friedrich’s positivity assumption or $\Omega$-filling advection, with which 
%With one of the following  two conditions, 
the norm equivalence \eqref{cond:norm-equivalence-lp-ar} also holds with $p=1$, see e.g. 
\cite{BocGun16,Guermond04lp}.  
\begin{itemize}
	\item  (Friedrich’s positivity assumption) There exists a constant $\mu_0>0$ such that 
	$\mu(x) - \frac{1}{p}\nabla \cdot \mathbf{b}\geq \mu_0$, a.e. $x$ in  $\Omega$. 
	\item ($\Omega$-filling advection). If $c=0$ and 
	$\nabla \cdot\mathbf{b}=0$, assume that the exists $z_{\pm}\in G^{\infty}_{\mathbf{b}}(\Omega)$ with $\norm{z_{\pm}}_{L^\infty}>0$ such that 
	$-\mathbf{b}\cdot \nabla z_{\pm} =p$ in $\Omega$ and 
	$z_{\pm}=0$ on $\partial\Omega_{\pm}$.%, where $p\in (1,\infty)$.
\end{itemize}

\begin{rem}[time-dependent advection-reaction equations]\label{rem:time-dependent-ar}
	The problem $\partial_t u + \mathbf{b}\cdot\nabla u + c u =f $
	with initial and boundary values can be recast into the form of 
	\eqref{eq:advection-reaction-inflow}.
	In fact, we may introduce the following notations $\tilde{x}=\left(x_{1}, x_{2}, \ldots, x_{n+1},t\right)$ and 
	$Q=\Omega \times (0, T),\, \partial Q=\partial \Omega \times (0, T)\cup (\Omega \times(\{T\} \cup\{0\})$. Also, we define the normal vector
	%\begin{equation*} 
	$\tilde{\mathbf{n}}= \left\{ \begin{array}{l } 
	(\mathbf{n}, 0) \text { on } \partial \Omega \times ( 0, T)\\
	(0,1) \text { on } \Omega \times\{T\} \\
	(0,-1) \text { on } \Omega \times\{0\}\\ 
	\end{array}\right..$
	% \end{equation*}
	Let $\tilde{\mathbf{b}}(\tilde{ {x}})=(\mathbf{b}(x , t), 1)$ and 
	$\partial Q_{-}=\{\tilde{\mathbf{x}} \in \partial Q \mid \tilde{\mathbf{b}} \cdot \tilde{\mathbf{n}} <0\}=
	(\partial \Omega_{-} \times [0, T]) \cup (\Omega \times\{0\})$, and 
	$\tilde{\nabla} \equiv\left(\nabla, \partial_{t}\right)$.  
	Then the problem can be written as 
	$\tilde{\mathbf{b}}\cdot \tilde{\nabla} u + c u = f$ with $u$ given on $\partial Q_{-}$. This reformulation has been used in many works, e.g.,  in \cite{PoiAze96}. %
\end{rem}

\iffalse 
With the norm equivalence \eqref{cond:norm-equivalence-lp-ar} and previous discussions on error estimates, we can readily establish  a prior and a posterior estimates for the following minimization problem
\begin{equation*}
    \inf_{v_{\bN,\theta}\in \mathfrak{N}_{\theta,n}}\mathcal{J}_\tau (v_{\bN,\theta}),
    \quad \mathcal{J}_\tau (v) = \norm{Av- f}_{L^p}^p + \tau \norm{Bv-g}_{L^p_{\abs{
    \mathbf{b}\cdot\mathbf{n}}}}^p.
\end{equation*}
As in Section \ref{sec:pinn-error}, we consider  the functional 
\begin{equation}
\mathcal{J}_\tau^\epsilon (v)= \norm{\varphi_\epsilon (\abs{Av -f})}_{L^1}  + \tau \norm{\varphi_{\epsilon_b} (\abs{Bv-g})\abs{\mathbf{b}\cdot \mathbf{n}}^{1/p}}_{L^1}.
\end{equation}
and a minimization problem of the functional over the space 
$X\cap \mathfrak{N}_{\theta,n}$, where  $\varphi_\epsilon (r)=r^{2m} (r+ \epsilon)^{p-2m}  $ for any $r>0$. 
Recall that $m\geq 1$ is an integer
 such that
 $p \in (2(m-1), 2m] \cap [1,\infty)$.
 \fi

%------------------------------------------------

\subsection{Integro-differential equations}\label{sec:ssec-fractional-laplacian}
Let $\Omega$ be a bounded Lipschitz domain satisfying the exterior ball condition or a $C^{2}$ bounded domain.  Consider the following  operator
%\begin{equation}
$    A   =   (- \Delta)^{\alpha/2}     +  \mathbf{b}\cdot \nabla   +    c  , \quad x\in \Omega \subset \Real^d, \quad 1<\alpha<2$ 
%\end{equation}
and $B =\text{Id} $ on $\Real^d\setminus\Omega$ and the image of $B$ has a compact support and the 
fractional Laplacian is defined  as a singular integral operator  on  $\Real^d$ (see e.g. 
\cite{Lischke2020fractionallap}) 
\begin{equation}\label{eq:fracLap-v1}  (- \Delta)^{\alpha/2}u(x) = c_{ d, \alpha} \int_{\mathbb{R}^d } \frac{u(x) - u(y)} {|x-y|^{d+\alpha}} \, d y,\quad c_{ d,\alpha} = \frac{2^\alpha\Gamma(\frac{\alpha+d}{2})}{\pi^{d/2} \abs{\Gamma(-\alpha/2)}
} .\end{equation}

Assume that there exists a constant $c_0>0$ such that %
$2c-\nabla\cdot \mathbf{b}\in L^\infty(\Omega)$ and  $2c-\nabla\cdot \mathbf{b}\geq  2c_0>0$.
By the Lax-Milgram Lemma, we can readily obtain the existence and uniqueness of a solution and there exists a constant $C_1>0$ that 
 \begin{equation}
   C_1 \norm{u}_{{H}^{\alpha/2}(\Real^d)} \leq \norm{A u}_{H^{-\alpha/2}(\Omega)} + \norm{Bu}_{H^{\alpha/2}(\Real^d\setminus\Omega)}.
 \end{equation} 
%Here $\widetilde{H}^s(\Omega) = \set{v\in H^s(\Real^d)| \text{ support of } v\subset\overline{\Omega}}$ with the norm $\norm{v}_{\widetilde{H}^s(\Omega)} =\norm{\tilde{v}}_{H^s(\Real^n)}$ where  $\tilde{v}=v $ on $\Omega$ and  is zero otherwise.  

The fractional Laplacian can be written as  %\begin{equation*} % (- \Delta)^{\alpha/2}u(x) =
 $c_{ d, \alpha} \int_{\mathbb{R}^d } \frac{2u(x) -u(x+y) - u(x-y)} {2|y|^{d+\alpha}} \, d y$.
%\end{equation*}
Then we obtain that when $\mathbf{b}$ and $c$ are in $L^\infty$, there exists a constant $C_2>0$ such that for $u\in C^2_c(\Real^d)$, 
%\begin{equation}
  $  \norm{A u}_{H^{-\alpha/2}(\Omega)} + \norm{Bu}_{H^{\alpha/2}(\Real^d\setminus\Omega)}
    \leq C_2\norm{u}_{C^2_c(\Real^d)}$.
%\end{equation} 
Here $C^2_c(\Real^d)$ is a subspace of $C^2(\Real^d)$ and elements in $C^2_c(\Real^d)$ are   compactly  supported.  

Assumption \ref{assum:norm-equivalence-generic} is verified with  $V={H}^{\alpha/2}(\Real^d)$, 
$Y=L^2(\Omega)$, $Z=H^{\alpha/2}(\Real^d\setminus\Omega)$, $X=C^2_c(\Real^d)$.

\begin{rem}When $Y=L^\infty(\Omega)$, norm relations in H\"older spaces  can be obtained using regularity results for the fractional Poisson equation  %$(-\Delta)^{\alpha/2}u =f$,
e.g.,  in  \cite{Gru15,RosSer14}.
\end{rem}